\documentclass[twoside,twocolumn,english,showkeys]{revtex4-2}
\usepackage{ae,aecompl}
\usepackage[T1]{fontenc}
\usepackage{geometry}
\usepackage{color}
\usepackage{babel}
\usepackage{amsmath}
\usepackage{amsthm}
\usepackage{amssymb}
\usepackage{graphicx}
\PassOptionsToPackage{normalem}{ulem}
\usepackage{ulem}
\usepackage[unicode=true,pdfusetitle,
 bookmarks=true,bookmarksnumbered=false,bookmarksopen=false,
 breaklinks=false,pdfborder={0 0 1},backref=false,colorlinks=false]
 {hyperref}

\makeatletter

\providecommand{\tabularnewline}{\\}

\numberwithin{equation}{section}
\numberwithin{figure}{section}
\theoremstyle{remark}
\newtheorem*{rem*}{\protect\remarkname}
\theoremstyle{plain}
\newtheorem{thm}{\protect\theoremname}
\theoremstyle{plain}
\newtheorem{cor}[thm]{\protect\corollaryname}

\@ifundefined{showcaptionsetup}{}{%
 \PassOptionsToPackage{caption=false}{subfig}}
\usepackage{subfig}
\makeatother

\providecommand{\corollaryname}{Corollary}
\providecommand{\remarkname}{Remark}
\providecommand{\theoremname}{Theorem}

\begin{document}
\title{Efficient relaxation scheme for the SIR and related compartmental
models}
\author{Vo Anh Khoa}
\affiliation{Department of Mathematics, Florida A\&M University, Tallahassee, FL
32307, USA}
\email{anhkhoa.vo@famu.edu, vakhoa.hcmus@gmail.com}

\author{Pham Minh Quan}
\affiliation{Department of Computer Science, Florida State University, Tallahassee,
FL 32306, USA}
\author{Ja'Niyah Allen}
\affiliation{Department of Mathematics, Florida A\&M University, Tallahassee, FL
32307, USA}
\author{Kbenesh W. Blayneh}
\affiliation{Department of Mathematics, Florida A\&M University, Tallahassee, FL
32307, USA}
\begin{abstract}
In this paper, we introduce a novel numerical approach for approximating
the SIR model in epidemiology. Our method enhances the existing linearization
procedure by incorporating a suitable relaxation term to tackle the
transcendental equation of nonlinear type. Developed within the continuous
framework, our relaxation method is explicit and easy to implement,
relying on a sequence of linear differential equations. This approach
yields accurate approximations in both discrete and analytical forms.
Through rigorous analysis, we prove that, with an appropriate choice
of the relaxation parameter, our numerical scheme is non-negativity-preserving
and globally strongly convergent towards the true solution. These
theoretical findings have not received sufficient attention in various
existing SIR solvers. We also extend the applicability of our relaxation
method to handle some variations of the traditional SIR model. Finally,
we present numerical examples using simulated data to demonstrate
the effectiveness of our proposed method.
\end{abstract}
\maketitle

\section{Introduction}

In recent years, the world has witnessed the devastating impact of
infectious diseases on a global scale. From the rapid spread of COVID-19
to the resurgence of long-standing ailments like measles and influenza,
understanding the dynamics of epidemics has become crucial for protecting
public health. To gain a deeper understanding of these intricate phenomena,
scientists have increasingly relied on mathematical modeling as an
influential tool for unraveling the complex mechanisms governing disease
transmission. Among the various models, the Susceptible-Infectious-Recovered
(SIR) model has emerged as a fundamental framework, providing valuable
insights into epidemic dynamics; cf. e.g. \citep{Earn2002,Khaleque2017,Kudryashov2021}
for its applications to modeling the influenza, Ebola and COVID-19.

The SIR model, initially proposed in the early 20th century by Kermack
and McKendrick \citep{1927}, has since been refined and adapted to
address contemporary challenges. This model effectively captures the
fundamental dynamics of epidemics by dividing a population into three
distinct compartments: susceptible individuals, infectious individuals,
and removed individuals. By considering the interactions between these
compartments, the SIR model takes into account various factors such
as transmission rates, removal rates, and the depletion of susceptible
individuals over time. The ``removals'' in this context encompass
individuals who are isolated, deceased, or have recovered and gained
immunity. Additionally, the model assumes that individuals who have
gained immunity or recovered enter a new category that is not susceptible
to the disease.

Consider a homogeneously mixed group of individuals of total size
$N\gg1$. Let $t\in\left(0,T\right)$ be the time variable with $T>0$
being the final time of observations. We take into account the following
functions:
\begin{align*}
S\left(t\right) & =\text{number of susceptibles at time }t,\\
I\left(t\right) & =\text{number of infectives at time }t,\\
R\left(t\right) & =\text{number of removals at time }t.
\end{align*}

Initiated, again, by Kermack and McKendrick in 1927 \citep{1927},
the evolutionary dynamics of these individuals can be modeled through
the following system of ordinary differential equations (ODEs):

\begin{equation}
\begin{cases}
I'\left(t\right)=\beta S\left(t\right)I\left(t\right)-\gamma I\left(t\right),\\
S'\left(t\right)=-\beta S\left(t\right)I\left(t\right),\\
R'\left(t\right)=\gamma I\left(t\right).
\end{cases}\label{eq:4}
\end{equation}
Here, the following assumptions are considered.

\textbf{(A1)} The total population size is always $N>0$, meaning
that $S\left(t\right)+I\left(t\right)+R\left(t\right)=N$ for all
$t$.

\textbf{(A2)} We know the infection rate $\beta>0$ from the infection
process, and the removal rate $\gamma>0$ from the removal process.

\textbf{(A3)} The initial conditions are $S\left(0\right)=n>1$, $I\left(0\right)=a\ge1$
and $R\left(0\right)=0$.

The explicit solution of the SIR model, despite its basic structure,
is widely known to be unattainable due to the exponential nonlinearity
of the transcendental equation governing removals. Consequently, numerous
numerical methods have been proposed to address this fundamental model.
The Taylor expansion method, initially employed by Kermack and McKendrick
in 1927, approximates the exponential term, leading to an approximate
analytic solution. This technique was utilized to simulate the plague
epidemic of 1905-1906 in Bombay, India, and has since become a tutorial
resource for students at both undergraduate and graduate levels; cf.
\citep{Caldwell2004}.

Over the years, many different methods have been studied to solve
the SIR model. Piyawong et al. \citep{Piyawong2003} introduced an
unconditionally convergent scheme that captures the long-term behavior
of the SIR model, offering improved numerical stability compared to
conventional explicit finite difference methods. Mickens \citep{Mickens2007}
and, recently, Conte et al. \citep{Conte2022} proposed and analyzed
stable nonstandard finite difference methods that effectively preserve
the positivity of the SIR solutions. Semi-analytical methods, such
as the Adomian decomposition approach \citep{Makinde2007}, have also
been proposed, alongside other methods cited therein, to derive approximate
analytical solutions. Additionally, the solutions of the SIR model
can be expressed in terms of the Lambert function, as demonstrated
in the publications by \citep{Barlow2020,Prodanov2021}. Furthermore,
an alternative approach involving parametric analysis has been employed
to obtain analytical solutions; see e.g. \citep{Harko2014}. While
most of these approaches are local approximations, a recent global
semi-analytical approach utilizing the Pad\'e approximation has been
presented in \citep{Chakir2023}.

While the above-mentioned approximation methods have demonstrated
numerical effectiveness, their convergence theories have received
limited investigation. Some publications discussing discrete methods
focus solely on stability analysis, leaving the global convergence
analysis unexplored. In this study, we propose a novel numerical approach
that guarantees global convergence. Our approach employs a relaxation
procedure, derived from the conventional linearization technique,
to approximate the SIR model in a continuous setting. Unlike the classical
version, our modified procedure introduces a relaxation term. In the
existing literature on partial differential equations (e.g., \citep{Khoa2020,Khoa2021,Mitra2019,Slodicka2001}
and references cited therein), this relaxation term mitigates the
local convergence issues encountered by conventional linearization
techniques such as Newton's method. Consequently, it permits to choose
an arbitrary starting point while guaranteeing the global convergence.
Within our specific context, the relaxation term facilitates capturing
the non-negativity of solutions while preserving global convergence.
In addition, by relying only on the dependence of the relaxation constant
on the removal rate, our approach accurately captures the long-term
behavior of the system. Besides, our explicit and easy-to-implement
approximate scheme is governed by a sequence of linear differential
equations. The desired approximate solution can be obtained discretely
or analytically based on individual preference.

Our paper is \textcolor{black}{four-fold}. In section \ref{sec:2},
we begin by revisiting the transcendental equation for removals and
discussing the essential properties of the SIR model. The latter part
of this section focuses on introducing the proposed relaxation scheme
and establishing its theoretical foundations. We prove that this scheme
is globally strongly convergent and preserves non-negativity. Additionally,
we derive an error estimate in $C^{0}$. In section \ref{sec:3},
we extend the applicability of our proposed scheme to some variants
of the SIR model. Then, to validate the effectiveness of our method,
numerical examples are presented in Section \ref{sec:4}. Finally,
\textcolor{black}{we provide some concluding remarks in section \ref{sec:45}},
and the appendix contains the proofs of our central theorems.

\section{Relaxation procedure\label{sec:2}}

\subsection{Transcendental equation revisited}

It is well known that the SIR model can be solved from a transcendental
differential equation. Here, we revisit how to get such an equation
to complete our analysis of the proposed scheme below. Let $\mu=\beta/\gamma$
be the reciprocal relative removal rate. By the second and third equations
of system (\ref{eq:4}), we have
\[
\frac{dS}{dR}=-\mu S,
\]
which leads to $\ln\left(S\right)=-\mu R+c$. Therefore, we arrive
at
\begin{equation}
S\left(t\right)=e^{c-\mu R\left(t\right)}.\label{eq:2}
\end{equation}
To find $c$, we set $t=0$ in (\ref{eq:2}) and use (A3). Indeed,
by $n=S\left(0\right)=e^{c-\mu R\left(0\right)}=e^{c}$, we get $c=\ln\left(n\right)$
and thus, deduce that
\begin{equation}
S\left(t\right)=ne^{-\mu R\left(t\right)}.\label{eq:3-1}
\end{equation}
Combining this and (A1), we derive the following nonlinear differential
equation for $R\left(t\right)$:
\begin{equation}
R'\left(t\right)=\gamma\left(N-ne^{-\mu R\left(t\right)}-R\left(t\right)\right).\label{eq:3}
\end{equation}

\begin{rem*}
To this end, the notation $C^{m}$ is used to denote the space of
functions with $m$ continuous derivatives. For the particular $C^{0}$
space, it is the space of continuous functions on $\left[0,T\right]$
with the standard max norm.
\end{rem*}
\begin{thm}
\label{thm:1}The differential equation (\ref{eq:3}) admits a unique
$C^{1}$ non-negative solution $R\left(t\right)$. Moreover, the existence
and uniqueness in $C^{1}$ of positive $S\left(t\right)$ and $I\left(t\right)$
to the SIR system (\ref{eq:4}) follow.
\end{thm}

\begin{proof}
The positivity of $S\left(t\right)$ and $I\left(t\right)$ is guaranteed
by the first and second equations of system (\ref{eq:4}), i.e.
\begin{align*}
S\left(t\right) & =n\exp\left(-\beta\int_{0}^{t}I\left(s\right)ds\right),\\
I\left(t\right) & =a\exp\left(\int_{0}^{t}\left(\beta S\left(s\right)-\gamma\right)ds\right).
\end{align*}
Then, by the third equation of (\ref{eq:4}), the non-negativity of
$R\left(t\right)$ follows. 

Cf. \citep[Theorem 3.2]{Sideris2013}, since the right hand side of
(\ref{eq:3}) is globally Lipschitzian, the equation admits a unique
local $C^{1}$ solution. Moreover, in view of the fact that $\left|\gamma\left(N-ne^{-\mu R\left(t\right)}-R\left(t\right)\right)\right|\le\gamma\left|R\left(t\right)\right|+\gamma\left(N+n\right)$
for any $t\ge0$, the obtained solution is global as a by-product
of \citep[Theorem 3.9]{Sideris2013}. Observe that the right hand
side of (\ref{eq:3-1}) is decreasing in the argument of $R\left(t\right)$.
Thereby, the existence and uniqueness of $S\left(t\right)$ follows.
We also get the existence and uniqueness of $I\left(t\right)$ in
view of the fact that the total population is conserved; cf. (A1).

Hence, we complete the proof of the theorem.
\end{proof}
Observe that if one can approximate $R\left(t\right)$ well in (\ref{eq:3}),
then $S\left(t\right)$ and $I\left(t\right)$ will be well approximated
via (\ref{eq:3-1}) and (\ref{eq:4}), respectively. Define $g\left(r\right)=\gamma ne^{-\mu r}+\gamma r$
for $r\in\mathbb{R}$. We can rewrite (\ref{eq:3}) as
\[
R'\left(t\right)=\gamma N-g\left(R\left(t\right)\right).
\]

\begin{rem*}
By the first and second equations of the SIR system (\ref{eq:4}),
we get
\[
\frac{dI}{dS}=-1+\frac{1}{\mu S}.
\]
Therefore, using $S\left(0\right)=n$ and $I\left(0\right)=a$, we
find that $I\left(t\right)-a=-S\left(t\right)+n+\frac{1}{\mu}\ln\left(S\left(t\right)\right)-\frac{1}{\mu}\ln\left(n\right)$.
Equivalently, we deduce that
\[
I\left(t\right)=\frac{1}{\mu}\ln\left(S\left(t\right)\right)-S\left(t\right)+a+n-\frac{1}{\mu}\ln\left(n\right).
\]
Since function $f\left(S\right)=\frac{1}{\mu}\ln\left(S\right)-S$
for $S>0$ attains its maximum at $S=\mu^{-1}$, we can estimate the
so-called amplitude, which is the maximum value of $I$, in the following
manner.
\begin{equation}
I_{\max}=-\frac{1}{\mu}\ln\left(\mu\right)-\frac{1}{\mu}+a+n-\frac{1}{\mu}\ln\left(n\right).\label{eq:I_max}
\end{equation}
\end{rem*}

\subsection{Derivation and analysis of the numerical scheme}

Let $\left\{ R_{k}\right\} _{k=0}^{\infty}$ be a time-dependent sequence
satisfying, for $k=1,2,3,\ldots$,
\begin{equation}
R_{k}'\left(t\right)+MR_{k}\left(t\right)=\gamma N-g\left(R_{k-1}\left(t\right)\right)+MR_{k-1}\left(t\right).\label{eq:5}
\end{equation}

The sequence $\left\{ R_{k}\right\} _{k=0}^{\infty}$ aims to approximate
$R\left(t\right)$ in (\ref{eq:3}) in the sense that $R_{k}$ will
be close to $R$ as $k\to\infty$ uniformly in time. The accompanying
initial condition for equation (\ref{eq:5}) is $R_{k}\left(0\right)=0$
for any $k\ge0$. Since our approximate model performs as an iterative
scheme, we need a starting point, $R_{0}\left(t\right)$. Here, we
choose $R_{0}\left(t\right)=0$ based on the initial condition of
$R\left(t\right)$ (cf. (A3)), which is the best information given
to the sought $R\left(t\right)$.

Also, in (\ref{eq:5}), we introduce a $k$-independent constant $M\ge\gamma>0$
for the so-called relaxation process. This relaxation term plays a
very important role. It allows us to prove the non-negativity of the
relaxation scheme, while many numerical approaches, including the
regular linearization method, do not have or cannot prove this feature.
Herewith, the regular linearization method we meant is the scheme
$\left\{ R_{k}\right\} _{k=0}^{\infty}$ in (\ref{eq:5}) with either
$M=0$ or only the exponential term in $g\left(r\right)$ being linearized.

Formulated below is the theorem showing that the scheme $\left\{ R_{k}\right\} _{k=0}^{\infty}$
preserves the non-negativity of the removals over the relaxation process.
\begin{thm}
\label{thm:2}The sequence $\left\{ R_{k}\right\} _{k=0}^{\infty}$
is a non-negativity-preserving scheme. Moreover, it holds true that
for all $M\ge\gamma$,
\[
0\le g\left(R_{k}\right)\le\gamma n+MR_{k}\quad\text{for any }k\ge0\;\text{and }t\ge0.
\]
\end{thm}

\begin{proof}
We prove this theorem by induction. The statement holds true for $k=1$.
Indeed, since $R_{0}\left(t\right)=0$, the equation for $R_{1}\left(t\right)$
reads as
\[
R_{1}'\left(t\right)+MR_{1}\left(t\right)=\gamma N-\gamma n\ge0.
\]
Therefore, we get
\begin{align*}
R_{1}\left(t\right) & =\frac{\gamma\left(N-n\right)}{M}\left(1-e^{-Mt}\right)\ge0,\\
g\left(R_{1}\right) & =\gamma ne^{-\mu R_{1}\left(t\right)}+\gamma R_{1}\left(t\right)\le\gamma n+MR_{1}\left(t\right).
\end{align*}
Next, assume that the statement holds true for $k=k_{0}$. We prove
that it also holds true for $k=k_{0}+1$. By (\ref{eq:5}), we have
\begin{align*}
R_{k_{0}+1}'\left(t\right) & +MR_{k_{0}+1}\left(t\right)\\
 & =\gamma N-g\left(R_{k_{0}}\left(t\right)\right)+MR_{k_{0}}\left(t\right)\ge0.
\end{align*}
Thus, we obtain $R_{k_{0}+1}\left(t\right)\ge e^{-Mt}R_{k_{0}+1}\left(0\right)\ge0.$
As a by product, we can estimate that
\begin{align*}
g\left(R_{k_{0}+1}\right) & =\gamma ne^{-\mu R_{k_{0}+1}\left(t\right)}+\gamma R_{k_{0}+1}\left(t\right)\\
 & \le\gamma n+MR_{k_{0}+1}\left(t\right).
\end{align*}
Hence, we complete the proof of the theorem.
\end{proof}
In the following, we formulate the strong convergence result for the
scheme $\left\{ R_{k}\right\} _{k=0}^{\infty}$. For ease of presentation,
proof of this result is deliberately placed in the Appendix. It is
worth mentioning that proof of the strong convergence of the scheme
relies so much on the strict estimation of $g'$. Such an estimation
can merely be obtained by the aid of the non-negativity of the scheme.
\begin{thm}
\label{thm:main}The sequence $\left\{ R_{k}\right\} _{k=0}^{\infty}$
defined in (\ref{eq:5}) is strongly convergent in $C^{0}$ toward
the true solution $R\left(t\right)$ to Equation (\ref{eq:3}). In
particular, we can find a number $C=C\left(T,M,\gamma,n,\mu\right)>0$
independent of $k$ such that the following error estimate holds true:
\[
\max_{0\le t\le T}\left|R_{k}\left(t\right)-R\left(t\right)\right|^{2}\le\frac{C^{k}}{k!}\max_{0\le t\le T}\left|R\left(t\right)\right|^{2}.
\]
\end{thm}

Our theoretical finding below shows that when $n\mu<1$, the scheme
$\left\{ R_{k}\left(t\right)\right\} _{k=0}^{\infty}$ converges faster
than the case $n\mu\ge1$. The proof of the following corollary is
also found in the Appendix.
\begin{cor}
\label{cor:4}Assume that $n\mu<1$. We can find a constant $c\in\left(0,1\right)$
independent of $k$ such that the following error estimate holds true:
\begin{equation}
\max_{0\le t\le T}\left|R_{k}\left(t\right)-R\left(t\right)\right|\le c^{k}\max_{0\le t\le T}\left|R\left(t\right)\right|.\label{eq:2.9}
\end{equation}
\end{cor}

As readily expected, for every step $k$, we obtain a non-homogeneous
differential equation that can be effectively approximated in the
discrete framework. Mimicking the proof of Theorem \ref{eq:3} and
using Theorem \ref{thm:2}, we have
\begin{equation}
R_{k}'\left(t\right)+MR_{k}\left(t\right)\le\gamma N-\gamma n+\left(M-\gamma\right)R_{k-1}\left(t\right),\label{eq:2.6}
\end{equation}
leading to
\begin{align*}
R_{k}\left(t\right) & \le e^{-Mt}\int_{0}^{t}e^{Ms}\left(\gamma a+\left(M-\gamma\right)R_{k-1}\left(s\right)\right)ds\\
 & \le t\gamma a+\left(M-\gamma\right)\int_{0}^{t}R_{k-1}\left(s\right)ds.
\end{align*}
By induction and by the choice $R_{0}\left(t\right)=0$, we can show
that
\[
\left|R_{k}\left(t\right)\right|\le\gamma a\sum_{i=1}^{k}\left(M-\gamma\right)^{i-1}\frac{t^{i}}{i!}.
\]
Therefore, if we choose $\gamma\le M\le\gamma+\frac{1}{T}$, then
$\left|R_{k}\left(t\right)\right|\le T\gamma a\left(e-1\right)$.
Note that this bound is independent of $k$. Thus, by Theorem \ref{thm:1},
we get $R_{k}\in C^{1}$ for any $k$ with, cf. (\ref{eq:2.6}), Theorem
\ref{thm:2} and the choice $M\le\gamma+\frac{1}{T}$,
\[
\left|R_{k}'\left(t\right)\right|\le\gamma a+\left(M-\gamma\right)T\gamma a\left(e-1\right)\le\gamma ae.
\]
Furthermore, by differentiating both sides of (\ref{eq:5}) with respect
to time, we can demonstrate that $R_{k}\in C^{2}$ for any $k$. Indeed,
\begin{align*}
R_{k}''\left(t\right) & +MR_{k}'\left(t\right)\\
 & =\left(M-\gamma+\gamma n\mu e^{-\mu R_{k-1}\left(t\right)}\right)R_{k-1}'\left(t\right)\\
 & \le\left(M-\gamma+\gamma n\mu\right)\left|R_{k-1}'\left(t\right)\right|.
\end{align*}
This yields that $\left|R_{k}''\left(t\right)\right|\le2M-\gamma+\gamma n\mu$,
which is a $k$-independent upper bound. This ensures the Euler method's
global error by leveraging its existing convergence theory. For completeness,
we present below the discrete solution to our proposed scheme.

Consider the time increment $\Delta t=T/P$ for $P\ge2$ being a fixed
integer. Then, we set the mesh-point in time by $t_{p}=p\Delta t$
for $0\le p\le P$. We seek $R_{k}^{p}\approx R_{k}\left(t_{p}\right)$
as a discrete solution to equation (\ref{eq:5}). By the standard
Euler method, $R_{k}^{p}$ is determined by the following equation:
\begin{align}
R_{k}^{p} & +\Delta tMR_{k}^{p}\nonumber \\
 & =R_{k}^{p-1}+\Delta t\left(\gamma N-g\left(R_{k-1}^{p}\right)+MR_{k-1}^{p}\right).\label{eq:dis}
\end{align}
By this way, the global error of the Euler method is attained in the
sense that for every $k$, there exists a constant $\tilde{C}>0$
such that
\[
\max_{0\le p\le P}\left|R_{k}^{p}-R_{k}\left(t_{p}\right)\right|\le\tilde{C}\Delta t.
\]
We accentuate that by the above analysis of $R_{k}$, i.e. $R_{k}\in C^{2}$
for any $k$, the constant $\tilde{C}$ is independent of $k$. Thus,
by Theorem \ref{thm:main}, we can estimate the distance between the
discrete (approximate) solution $R_{k}^{p}$ and the true solution
$R$ at each mesh-point, 
\begin{equation}
\max_{0\le p\le P}\left|R_{k}^{p}-R\left(t_{p}\right)\right|\le\tilde{C}\Delta t+\sqrt{\frac{C^{k}}{k!}}\max_{0\le t\le T}\left|R\left(t\right)\right|.\label{eq:2.7-1}
\end{equation}

\begin{rem*}
We have the following remarks:
\begin{itemize}
\item After obtaining the approximator $R_{k}^{p}$ for $R\left(t_{p}\right)$,
we can compute $S\left(t_{p}\right)$ using (\ref{eq:3-1}). Then,
the approximate solution for $I\left(t_{p}\right)$ can be determined
using (A1), specifically $I\left(t_{p}\right)=N-S\left(t_{p}\right)-R\left(t_{p}\right)$.
\item Both $\tilde{C}$ and $C$ in (\ref{eq:2.7-1}) are independent of
$P$ and $k$. As a by product, our discrete relaxation scheme $\left\{ R_{k}^{p}\right\} _{k=0}^{\infty}$,
as defined in (\ref{eq:dis}), is globally strongly convergent in
$C^{0}$. Similar to the proof of Theorem \ref{thm:2}, we can demonstrate
that $R_{k}^{p}$ is non-negativity-preserving. It is important to
note that many previous approximations, such as the method of series
expansions \citep{Makinde2007,Chakir2023}, parametrization method
\citep{Harko2014,Prodanov2023} and finite difference method \citep{Piyawong2003,Side2018},
did not adequately address the preservation of non-negativity/positivity.
Furthermore, certain recent positive numerical schemes fail to provide
an error bound, as observed in the publications \citep{Khalsaraei2021,Conte2022}.
\item While the application of the Euler method is initiated, we can further
employ higher-order numerical methods to produce a faster convergent
solver for the linear differential equation of $R_{k}$. Among these,
the Runge-Kutta method stands out as the most favorable choice, offering
a convergence rate of order $q\ge2$. Building upon the analysis of
the Euler method above, we can prove that all derivatives of the right
hand side of (\ref{eq:5}) exist up to order $q$ and $R_{k}\in C^{q}$
for any $k$. Therefore, we can show that $R_{k}^{p}$ globally converges
to $R_{k}$ with a rate of $\mathcal{O}\left(\Delta t^{q}\right)$;
cf. e.g. \citep[Theorem 3.4]{1993} for the existing theory on the
global convergence of the Runge-Kutta method. Note here that $\mathcal{O}\left(x\right)$
is the conventional Landau symbol.
\item Similar to (\ref{eq:2.7-1}), the convergence of the Runge-Kutta method
remains unaffected by $k$. Nevertheless, it is important to emphasize
that this convergence is heavily contingent upon the upper bounds
of the involved derivatives. Considering the boundedness of $R_{k},R_{k}'$
and $R_{k}''$ established above, it becomes evident that these bounds
tend to increase as the order rises. Consequently, it is crucial to
exclusively employ variants of the Runge-Kutta method with appropriately
high orders. This perspective holds true when applying the Runge-Kutta
method directly to the differential equation (\ref{eq:3}).
\end{itemize}
\end{rem*}

\section{Extensions to other SIR models\label{sec:3}}

In this section, we briefly discuss the applicability of the relaxation
method to other population models of SIR type. In particular, we show
below how the proposed approach can be adapted to approximate the
SIRD (Susceptible-Infectious-Recovered-Deceased) and SIRX models;
cf. \citep{Bailey1975,Maier2020,Baerwolff2021} for an overview of
these models.

\subsection*{SIRD model}

The SIRD model extends the SIR model by distinguishing between recovered
and deceased individuals. In this framework, the removals in the SIR
model no longer encompass the number of infected individuals who have
passed away. To account for mortality, a mortality rate $\sigma>0$
is introduced, representing the rate at which infected individuals
succumb to death. Consequently, the death rate per unit of time is
calculated as the product of the mortality rate and the number of
infected individuals. Additionally, as the number of deceased individuals
is excluded from the removals, the rate of change of infections over
time is adjusted to reflect the loss caused by mortality. Mathematically,
the SIRD model reads as
\begin{equation}
\begin{cases}
I'\left(t\right)=\beta S\left(t\right)I\left(t\right)-\left(\gamma+\sigma\right)I\left(t\right),\\
S'\left(t\right)=-\beta S\left(t\right)I\left(t\right),\\
R'\left(t\right)=\gamma I\left(t\right),\\
D'\left(t\right)=\sigma I\left(t\right),
\end{cases}\label{eq:3.1}
\end{equation}
where $D\left(t\right)$ stands for the number of deceased people
(after infection) at time $t$. In (\ref{eq:3.1}), we make use of
the following assumptions.

\textbf{(B1)} The total population size is always conserved with $N>0$,
meaning that $S\left(t\right)+I\left(t\right)+R\left(t\right)+D\left(t\right)=N$
for all $t$.

\textbf{(B2)} We know the infection rate $\beta>0$ from the infection
process, the removal rate $\gamma>0$ from the removal process, and
the death rate $\sigma>0$ from the mortality process.

\textbf{(B3)} The initial conditions are $S\left(0\right)=n>1$, $I\left(0\right)=a\ge1$,
$R\left(0\right)=0$ and $D\left(0\right)=0$.
\begin{rem*}
By the first and second equations of the SIRD system (\ref{eq:3.1}),
we see that
\[
\frac{dI}{dS}=-1+\frac{\gamma+\sigma}{\beta S}.
\]
Similar to the classic SIR model (\ref{eq:4}), we can thus formulate
the so-called amplitude in the following fashion:
\begin{align}
I_{\max} & =\frac{\gamma+\sigma}{\beta}\ln\left(\frac{\gamma+\sigma}{\beta}\right)\nonumber \\
 & -\frac{\gamma+\sigma}{\beta}+a+n-\frac{\gamma+\sigma}{\beta}\ln\left(n\right)\label{eq:Imax-SIRD}
\end{align}
when $S$ reaches $\frac{\gamma+\sigma}{\beta}$.
\end{rem*}
\begin{rem*}
The SIRD model (\ref{eq:3.1}) resembles the SIRX model without containment
rate. In the SIRX model, an additional class called ``$X$\textquotedbl{}
was introduced to account for the impact of social or individual behavioral
changes during quarantine. Individuals in this class, referred to
as symptomatic quarantined individuals, no longer contribute to the
transmission of the infection. Instead of $\sigma$, the SIRX model
without containment rate consideres $\kappa>0$ that represents the
rate at which infected individuals are removed due to quarantine measures.
The SIRX model with the containment rate is not the scope of our paper
since the associated transcendental system does not take the same
form of (\ref{eq:5}). Indeed, the transcendental system governing
the full SIRX model is of an integro-differential equation.
\end{rem*}
Now, we detail the transcendental equation for $R\left(t\right)$
and the application of the relaxation scheme. From the second and
third equations of system (\ref{eq:3.1}), we deduce that
\begin{equation}
S\left(t\right)=ne^{-\mu R\left(t\right)},\label{eq:3.2}
\end{equation}
where we have recalled the reciprocal relative removal constant $\mu=\beta/\gamma$.
Using the same way, the third and last equations of system (\ref{eq:3.1})
give
\begin{equation}
D\left(t\right)=\frac{\sigma}{\gamma}R\left(t\right)\label{eq:3.3}
\end{equation}
by virtue of $R\left(0\right)=D\left(0\right)=0$ (cf. (B3)). Then,
plugging (\ref{eq:3.2}), (\ref{eq:3.3}) and (B1) into the third
equation of (\ref{eq:3.1}), we obtain the following differential
equation for $R\left(t\right)$:
\begin{equation}
R'\left(t\right)=\gamma\left[N-ne^{-\mu R\left(t\right)}-\left(1+\frac{\sigma}{\gamma}\right)R\left(t\right)\right].\label{eq:3.4}
\end{equation}
Henceforth, our relaxation scheme in this case becomes
\begin{equation}
R_{k}'\left(t\right)+\overline{M}R_{k}\left(t\right)=\gamma N-\overline{g}\left(R_{k-1}\left(t\right)\right)+\overline{M}R_{k-1}\left(t\right),\label{eq:3.5}
\end{equation}
where $\overline{g}\left(r\right)=\gamma ne^{-\mu r}+\left(\gamma+\sigma\right)r$.
Similar to the SIR model, here we rely on (B3) to choose $R_{k}\left(0\right)=0$
for any $k\ge0$ as the initial condition and $R_{0}\left(t\right)=0$
as the starting point.

By choosing $\overline{M}\ge\gamma+\sigma$, our sequence $\left\{ R_{k}\right\} _{k=0}^{\infty}$
(defined in (\ref{eq:3.5})) is non-negativity-preserving and globally
strongly convergent to $R$ of the transcendental equation (\ref{eq:3.4}).
These findings are analogous to our central Theorems \ref{thm:2}
and \ref{thm:main}, and therefore, we omit their formulations. Besides,
Theorem \ref{thm:1} is applied to (\ref{eq:3.4}), guaranteeing the
global existence and uniqueness of the $C^{1}$ solutions to the SIRD
model (\ref{eq:3.1}).

\subsection*{SIR model with background mortality}

The SIR model, along with its variants SIRD and SIRX, assumes a constant
population size. These models, known as epidemiological SIR-type models
without vital dynamics, are limited in their representation of population
changes; see (A1), (B1) and (C1). The SIR model with vital dynamics
addresses this limitation by incorporating birth and death rates to
account for population size fluctuations.

In the present work, we explore that the transcendental system governing
the SIR model with background mortality takes the form of (\ref{eq:5}).
With $\sigma$ being the death rate, the population experiences changes
over time. Here, individuals from all compartments can exit through
deaths, allowing for a more realistic representation of population
dynamics. Mathematically, the SIR model with background mortality
can be expressed as follows:

\begin{equation}
\begin{cases}
I'\left(t\right)=\beta S\left(t\right)I\left(t\right)-\gamma I\left(t\right)-\sigma I\left(t\right),\\
S'\left(t\right)=-\beta S\left(t\right)I\left(t\right)-\sigma S\left(t\right),\\
R'\left(t\right)=\gamma I\left(t\right)-\sigma R\left(t\right).
\end{cases}\label{eq:3.13-1}
\end{equation}

In this perspective, we make use of the following assumptions.

\textbf{(C1)} The total population size is dependent of $t$, i.e.
$N=N\left(t\right)>0$. It can be computed that $N\left(t\right)=S\left(t\right)+I\left(t\right)+R\left(t\right)=e^{-\sigma t}N_{0}$
for some fixed $N_{0}>0$.

\textbf{(C2)} We know the infection rate $\beta>0$ from the infection
process, the removal rate $\gamma>0$ from the removal process, and
the death rate $\sigma>0$ from the mortality process.

\textbf{(C3)} The initial conditions are $S\left(0\right)=n>1$, $I\left(0\right)=a\ge1$
and $R\left(0\right)=0$. This implies that $N_{0}=n+a$.

Similar to the classical SIR model, we seek the transcendental equation
for $R$ prior to the application of our proposed numerical scheme.
When doing so, we define $\overline{R}\left(t\right)=e^{\sigma t}R\left(t\right)$
and $\overline{S}\left(t\right)=e^{\sigma t}S\left(t\right)$. By
the second and third equations of system (\ref{eq:3.13-1}), we find
that
\begin{align}
\overline{S}'\left(t\right) & =-\beta\overline{S}\left(t\right)I\left(t\right),\label{eq:3.14}\\
\overline{R}'\left(t\right) & =\gamma e^{\sigma t}I\left(t\right).\label{eq:3.15}
\end{align}
Therefore, we deduce that
\begin{align}
\frac{d\overline{S}}{d\overline{R}} & =-\mu e^{-\sigma t}\overline{S},\label{eq:3.181}
\end{align}
or equivalently, $\ln\left(\overline{S}\right)=-\mu e^{-\sigma t}\overline{R}+\tilde{c}\left(t\right)$.
Herewith, we have recalled the reciprocal relative removal rate $\mu=\beta/\gamma$.
Henceforth, we have
\begin{equation}
\overline{S}\left(t\right)=e^{\tilde{c}\left(t\right)-\mu e^{-\sigma t}\overline{R}\left(t\right)}.\label{eq:3.18-1}
\end{equation}
Since $\overline{S}\left(0\right)=n$ and $\overline{R}\left(0\right)=R\left(0\right)=0$
by (D3), we find that $\tilde{c}\left(0\right)=\ln\left(n\right)$.
Moreover, by taking the derivative in time of (\ref{eq:3.18-1}),
we arrive at
\begin{equation}
\overline{S}'\left(t\right)=e^{\tilde{c}\left(t\right)-\mu e^{-\sigma t}\overline{R}\left(t\right)}\left[\tilde{c}'\left(t\right)-\mu e^{-\sigma t}\left(-\sigma+\overline{R}'\left(t\right)\right)\right].\label{eq:3.19}
\end{equation}
Dividing both sides of (\ref{eq:3.19}) by $\overline{R}'\left(t\right)$,
we find that

\begin{align*}
 & \frac{\overline{S}'\left(t\right)}{\overline{R}'\left(t\right)}\\
 & =\frac{e^{\tilde{c}\left(t\right)-\mu e^{-\sigma t}\overline{R}\left(t\right)}\left[\tilde{c}'\left(t\right)-\mu e^{-\sigma t}\left(-\sigma+\overline{R}'\left(t\right)\right)\right]}{\overline{R}'\left(t\right)}.
\end{align*}
Then combining this with (\ref{eq:3.181}) and (\ref{eq:3.18-1}),
we derive the following differential equation for $\tilde{c}\left(t\right)$:
\begin{align*}
 & e^{\tilde{c}\left(t\right)-\mu e^{-\sigma t}\overline{R}\left(t\right)}\left[\tilde{c}'\left(t\right)-\mu e^{-\sigma t}\left(-\sigma+\overline{R}'\left(t\right)\right)\right]\\
 & =-\mu e^{-\sigma t}e^{\tilde{c}\left(t\right)-\mu e^{-\sigma t}\overline{R}\left(t\right)}\overline{R}'\left(t\right),
\end{align*}
or equivalently, $\tilde{c}'\left(t\right)=-\mu\sigma e^{-\sigma t}$.
Thus, we obtain $\tilde{c}\left(t\right)=\ln\left(n\right)+\mu\left(e^{-\sigma t}-1\right)$
and
\[
\overline{S}\left(t\right)=ne^{\mu\left(e^{-\sigma t}-1\right)}e^{-\mu e^{-\sigma t}\overline{R}\left(t\right)}.
\]
Together with the back-substitution $e^{-\sigma t}\overline{R}\left(t\right)=R\left(t\right)$,
we thereby get $\overline{S}\left(t\right)=ne^{\mu\left(e^{-\sigma t}-1\right)}e^{-\mu R\left(t\right)}$.
Now, we note that by (C1) and (C3), $e^{\sigma t}I\left(t\right)=N_{0}-\overline{S}\left(t\right)-\overline{R}\left(t\right)$
holds true for any $t$. Plugging this into (\ref{eq:3.15}) and using
the fact that $\overline{S}\left(t\right)=ne^{\mu\left(e^{-\sigma t}-1\right)}e^{-\mu R\left(t\right)}=ne^{\mu\left(e^{-\sigma t}-1\right)}e^{-\mu e^{-\sigma t}\overline{R}\left(t\right)}$,
we derive the transcendental equation for $\overline{R}$ as follows:

\begin{equation}
\overline{R}'\left(t\right)=\gamma\left[N_{0}-ne^{\mu\left(e^{-\sigma t}-1\right)}e^{-\mu e^{-\sigma t}\overline{R}\left(t\right)}-\overline{R}\left(t\right)\right].\label{eq:3.18}
\end{equation}
By setting $\hat{g}\left(r\right)=\gamma ne^{\mu\left(e^{-\sigma t}-1\right)}e^{-\mu e^{-\sigma t}r}+\gamma r$,
our relaxation scheme for the SIR model with background mortality
is structured by
\begin{align}
\overline{R}_{k}'\left(t\right) & +\hat{M}\overline{R}_{k}\left(t\right)\nonumber \\
 & =\gamma N_{0}-\hat{g}\left(\overline{R}_{k-1}\left(t\right)\right)+\hat{M}\overline{R}_{k-1}\left(t\right).\label{eq:3.13-2}
\end{align}

Similar to the above-mentioned SIR-based models, we choose $\overline{R}_{k}\left(0\right)=0$
for any $k\ge0$ as the initial condition and $\overline{R}_{0}\left(t\right)=0$
as the starting point, based on the fact that $\overline{R}\left(0\right)=R\left(0\right)=0$.
Also, here we take $\hat{M}\ge\gamma$ to ensure the non-negativity
preservation and global strong convergence of the sequence $\left\{ \overline{R}_{k}\right\} _{k=0}^{\infty}$
(defined in (\ref{eq:3.13-2})) to the sought $\overline{R}$ of the
transcendental equation (\ref{eq:3.18}). As another analog of Theorems
\ref{thm:2} and \ref{thm:main}, we omit details of the formulations
of the theoretical results for the sequence $\left\{ \overline{R}_{k}\right\} _{k=0}^{\infty}$.
It is also worth mentioning that Theorem \ref{thm:1} remains true
in this case, providing the global existence and uniqueness of $C^{1}$
solutions to the SIR model with mortality (\ref{eq:3.13-1}).

\section{Numerical experiments\label{sec:4}}

In this section, we verify the numerical performance of the proposed
relaxation method. Initially, we employ various approaches to solve
the SIR model (\ref{eq:4}) for the purpose of comparison. These include
our method (\ref{eq:5}), as well as the standard methods: approximate
analytic solution, regular linearization procedure, and conventional
explicit Euler method. It is important to note that since the conventional
explicit Euler method is considered in this comparison, we also apply
the Euler method to our relaxation scheme, as outlined in (\ref{eq:dis}),
as well as to the regular linearization procedure.

Additionally, it is worth mentioning that the approximate analytic
solution for $R\left(t\right)$ (referred to as $R_{\text{a}}$) can
be found in \citep{1927,Caldwell2004,Baerwolff2021}. In particular,
it is of the following form: 
\begin{equation}
R_{\text{a}}\left(t\right)=\frac{1}{n\mu^{2}}\left[n\mu-1+\eta\tanh\left(\frac{\gamma\eta t}{2}-\psi\right)\right],\label{eq:4.1}
\end{equation}
where
\begin{align*}
\eta & =\left[2n\mu^{2}\left(N-n\right)+\left(n\mu-1\right)^{2}\right]^{1/2},\\
\psi & =\tanh^{-1}\left[\frac{1}{\eta}\left(n\mu-1\right)\right].
\end{align*}

The approximate analytic solution mentioned above corresponds to the
solution of the Riccati equation. However, it is applicable only when
$\mu R$ is sufficiently small. Furthermore, the conventional explicit
Euler method is expressed as follows:

\begin{equation}
R^{p}=R^{p-1}+\Delta t\left(\gamma N-g\left(R^{p-1}\right)\right),\label{eq:4.2}
\end{equation}
where $R^{p}\approx R\left(t_{p}\right)$ is the discrete solution
to the nonlinear differential equation (\ref{eq:3}) for $t_{p}=p\Delta t$
being the mesh-point in time.

In the second test, we utilize the widely used Runge-Kutta RK4 method
to solve the SIR model (\ref{eq:5}) by applying it to our relaxation
scheme (\ref{eq:4}). We then compare its performance when using the
Euler-relaxation method (\ref{eq:dis}) and when directly applying
the Runge-Kutta RK4 method to (\ref{eq:3}). For sake of clarity,
we provide the formulation of the RK4 method for solving a differential
equation of a general type $R'\left(t\right)=F\left(t,R\left(t\right)\right)$:

\begin{align}
R^{p} & =R^{p-1}+\frac{1}{6}K_{1}\left(t_{p-1},R^{p-1}\right)+\frac{1}{3}K_{2}\left(t_{p-1},R^{p-1}\right)\nonumber \\
 & +\frac{1}{3}K_{3}\left(t_{p-1},R^{p-1}\right)+\frac{1}{6}K_{4}\left(t_{p-1},R^{p-1}\right),\label{eq:RK4}
\end{align}
where we have denoted the intermediate values by
\begin{align}
K_{1}\left(t_{p-1},R^{p-1}\right) & =\Delta tF\left(t_{p-1},R^{p-1}\right),\label{eq:K}\\
K_{2}\left(t_{p-1},R^{p-1}\right) & =\Delta tF\left(t_{p-1}+\frac{\Delta t}{2},R^{p-1}+\frac{K_{1}}{2}\right),\\
K_{3}\left(t_{p-1},R^{p-1}\right) & =\Delta tF\left(t_{p-1}+\frac{\Delta t}{2},R^{p-1}+\frac{K_{2}}{2}\right),\\
K_{4}\left(t_{p-1},R^{p-1}\right) & =\Delta tF\left(t_{p-1}+\Delta t,R^{p-1}+K_{3}\right).\label{eq:KK}
\end{align}

\begin{widetext}
When using our method, it is important to note that for each iteration
$k$, we solve the linear differential equation $R_{k}'\left(t\right)=F\left(t,R_{k}\left(t\right)\right)$,
where $F\left(t,R_{k}\left(t\right)\right)=-MR_{k}\left(t\right)+\gamma N-g\left(R_{k-1}\left(t\right)\right)+MR_{k-1}\left(t\right)$.
Notice that in this equation, the presence of the midpoint $t_{p-1}+\frac{\Delta t}{2}$,
applied to $R_{k-1}\left(t\right)$ obtained from the previous step,
leads to the following linear approximation:

\begin{equation}
R_{k-1}\left(t_{p-1}+\frac{\Delta t}{2}\right)=\frac{1}{2}\left[R_{k-1}\left(t_{p-1}\right)+R_{k-1}\left(t_{p-1}+\Delta t\right)\right].\label{eq:RR}
\end{equation}
Denote this approximation by $R_{k-1}^{p-0.5}=R_{k-1}\left(t_{p-1}+\frac{\Delta t}{2}\right)$.
Thereby, we seek $R_{k}^{p}$ satisfying (\ref{eq:RK4}) in which
the intermediate values are given by
\begin{align}
K_{1} & =\Delta t\left(-MR_{k}^{p-1}+\gamma N-g\left(R_{k-1}^{p-1}\right)+MR_{k-1}^{p-1}\right),\\
K_{2} & =\Delta t\left[-M\left(R_{k}^{p-1}+\frac{K_{1}}{2}\right)+\gamma N-g\left(R_{k-1}^{p-0.5}\right)+MR_{k-1}^{p-0.5}\right],\\
K_{3} & =\Delta t\left[-M\left(R_{k}^{p-1}+\frac{K_{2}}{2}\right)+\gamma N-g\left(R_{k-1}^{p-0.5}\right)+MR_{k-1}^{p-0.5}\right],\\
K_{4} & =\Delta t\left[-M\left(R_{k}^{p-1}+K_{3}\right)+\gamma N-g\left(R_{k-1}^{p}\right)+MR_{k-1}^{p}\right].\label{eq:RRR}
\end{align}
\end{widetext}

On the other hand, when directly applying the Runge-Kutta RK4 method
to the nonlinear differential equation (\ref{eq:3}), we have $F\left(t,R\left(t\right)\right)=\gamma\left(N-ne^{-\mu R\left(t\right)}-R\left(t\right)\right)$.

In the third test, we present the numerical performance of the proposed
method in solving the SIR-based models discussed in section \ref{sec:3}.
In particular, our focus in this test is on
\begin{enumerate}
\item the scheme $\left\{ R_{k}\right\} _{k=0}^{\infty}$, defined in (\ref{eq:3.5}),
for the SIRD model (\ref{eq:3.1}). In this model, the relaxation
parameter satisfies $\overline{M}\ge\gamma+\sigma$.
\item the scheme $\left\{ R_{k}\right\} _{k=0}^{\infty}$ computed from
$\left\{ \overline{R}_{k}\right\} _{k=0}^{\infty}$ (defined in (\ref{eq:3.13-2}))
for the SIR model with background mortality (\ref{eq:3.13-1}). In
this case, we condition that $\hat{M}\ge\gamma$.
\end{enumerate}
To evaluate the accuracy of the relaxation scheme, we assess the proximity
of the approximation when approaching the maximum value of $I$. It
is important to recall that explicit expressions for $I_{\max}$ have
been derived for each specific case. The reader is referred to (\ref{eq:I_max})
for the classical SIRD model, and (\ref{eq:Imax-SIRD}) for the SIRD
model. For the SIR model with background mortality, since the maximum
value of $I$ cannot be found explicitly, we run the simulation with
several values of $P$ and $K$ to verify the numerical stability.
When increasing these parameters, we also identify the numerical amplitude
and peak day to see the performance of our relaxation method in the
Euler and RK4 frameworks.

\subsection*{Test 1}

In this test, we compare our Euler-relaxation approach with the approximate
analytic solution (\ref{eq:4.1}), the regular linearization procedure
(which arises when the relaxation parameter vanishes), and the direct
explicit Euler method (\ref{eq:4.2}). Alongside assessing numerical
stability, we evaluate the performance of these methods based on the
amplitude $I_{\max}$ presented in (\ref{eq:I_max}) and the peak
day.

Here, we consider a population sample of $N=1000$ for the SIR model
(\ref{eq:4}) over the course of one year ($T=365$). Initially, we
assume that there are $a=2$ infected people in this sample, leaving
$n=998$ individuals susceptible to infection. Furthermore, we choose
a removal rate of $\gamma=0.02$ and an infection rate of $\beta=0.0004$.
With these choices, we obtain a reciprocal relative removal rate of
$\mu=\beta/\gamma=0.02$, indicating that $n\mu=19.96>1$. Additionally,
for our relaxation process, we set $M=0.02$.

\begin{table}[t]
\begin{centering}
\begin{tabular}{|c|c|c|c|c|c|c|c|}
\hline 
Method \# & 1 & 1 & 2 & 2 & 3 & 4 & 4\tabularnewline
\hline 
\# of time step $P$ & $100$ & $1000$ & $100$ & $1000$ & None & $100$ & $1000$\tabularnewline
\hline 
\# of iteration $K$ & 5 & 50 & 2 & 4 & None & None & None\tabularnewline
\hline 
Amplitude $I_{\max}$ & 797 & 800 & 800 & 800 & 755 & 793 & 800\tabularnewline
\hline 
Peak day & 25 & 23 & 138 & 35 & 114 & 32 & 25\tabularnewline
\hline 
\end{tabular}
\par\end{centering}
\caption{Values of the computed amplitude $I_{\max}$ obtained from different
methods and the corresponding peak days. Method \#1: our Euler-relaxation
method (\ref{eq:dis}). Method \#2: the regular linearization method,
i.e. our proposed method (\ref{eq:5}) but with $M=0$. Method \#3:
the approximate analytic solution $R_{\text{a}}$ formulated in (\ref{eq:4.1}).
Method \#4: the conventional explicit Euler method (\ref{eq:4.2})
applied directly to the nonlinear differential equation (\ref{eq:3}).
By (\ref{eq:I_max}), the true amplitude $I_{\max,\text{true}}$ is
800 in this scenario.\label{tab:1}}
\end{table}

Our numerical results for Test 1 are presented in Table \ref{tab:1}.
Based on the maximum amplitude ($I_{\max}$), our proposed method
within the Euler context (method \#1) outperforms the approximations
obtained from methods \#2\textendash 4. The first two columns of Table
\ref{tab:1} demonstrate the numerical stability of our proposed method,
particularly when dealing with relatively small values of $P$ and
$K$. Remarkably, when $P=100$ and $K=5$, our method yields an $I_{\max}$
value of 797, which is very close to the true value of 800 as shown
in (\ref{eq:I_max}). In contrast, the $I_{\max}$ obtained from the
approximate analytic solution (method \#3) shows a significant deviation.
We also observe that the amplitude $I_{\max}$ obtained from method
\#3 remains unaffected regardless of the choice of $P$.

\begin{widetext}
\begin{figure*}
\begin{centering}
\subfloat[Method \#1 ($P=100,K=5$)]{\begin{centering}
\includegraphics[scale=0.18]{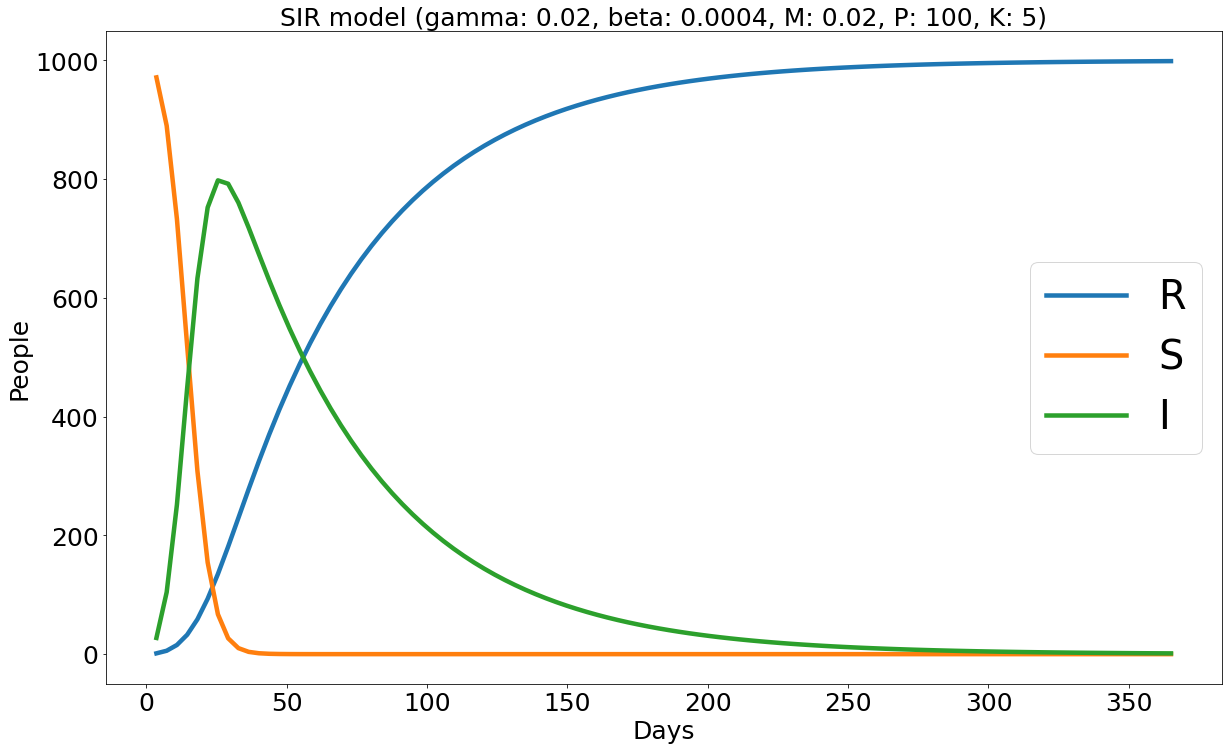}
\par\end{centering}

}\subfloat[Method \#1 ($P=1000,K=50$)]{\begin{centering}
\includegraphics[scale=0.18]{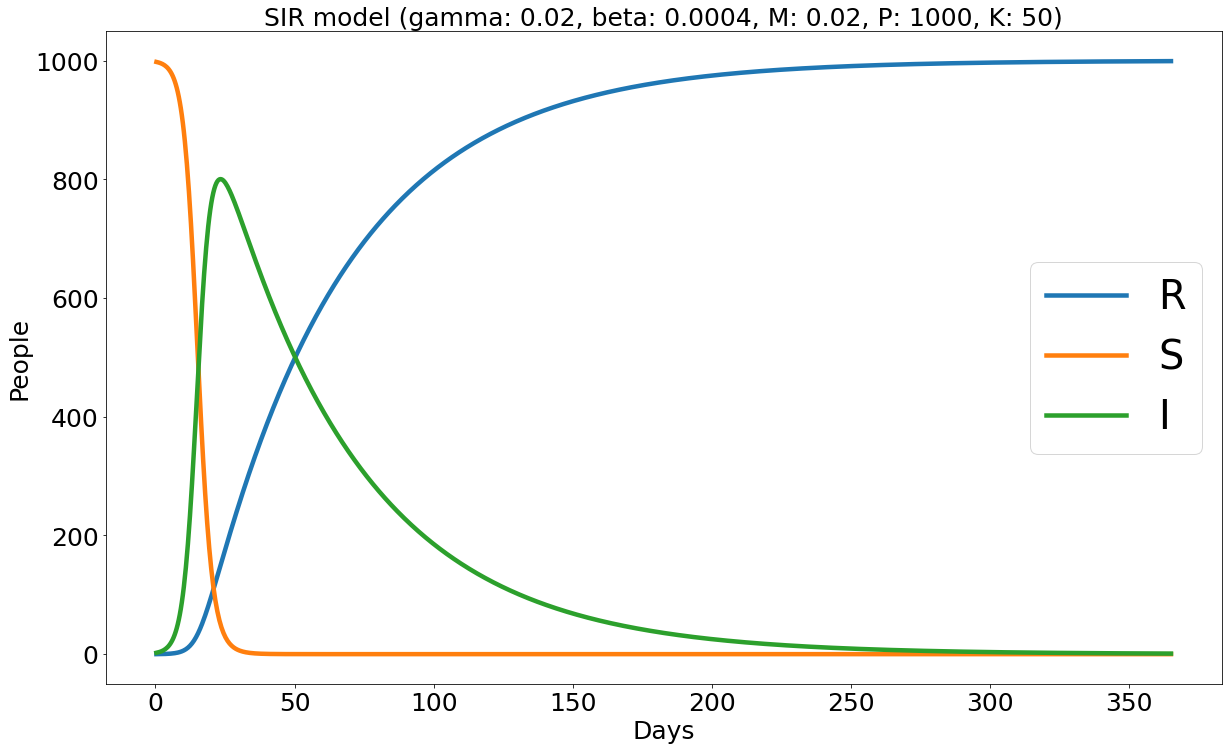}
\par\end{centering}
}
\par\end{centering}
\begin{centering}
\subfloat[Method \#2 ($P=100,K=2$)]{\begin{centering}
\includegraphics[scale=0.18]{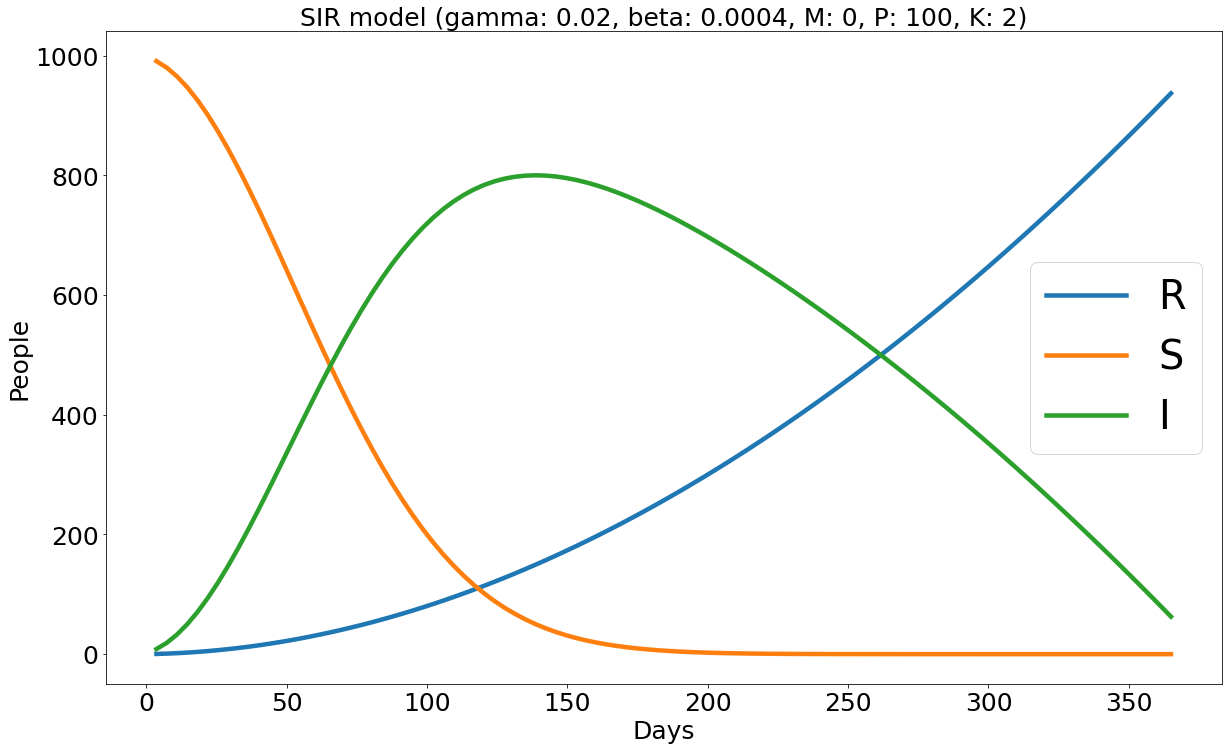}
\par\end{centering}
}\subfloat[Method \#2 ($P=1000,K=4$)]{\begin{centering}
\includegraphics[scale=0.18]{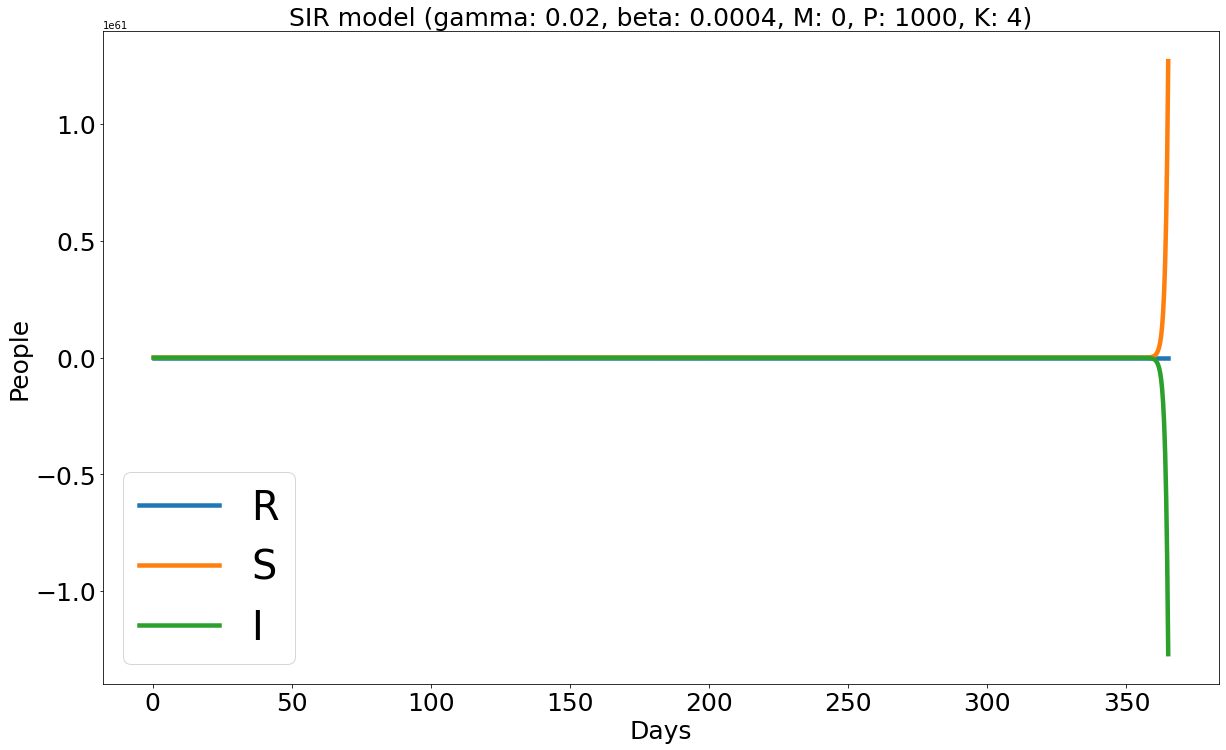}
\par\end{centering}
}
\par\end{centering}
\begin{centering}
\subfloat[Method \#3]{\begin{centering}
\includegraphics[scale=0.18]{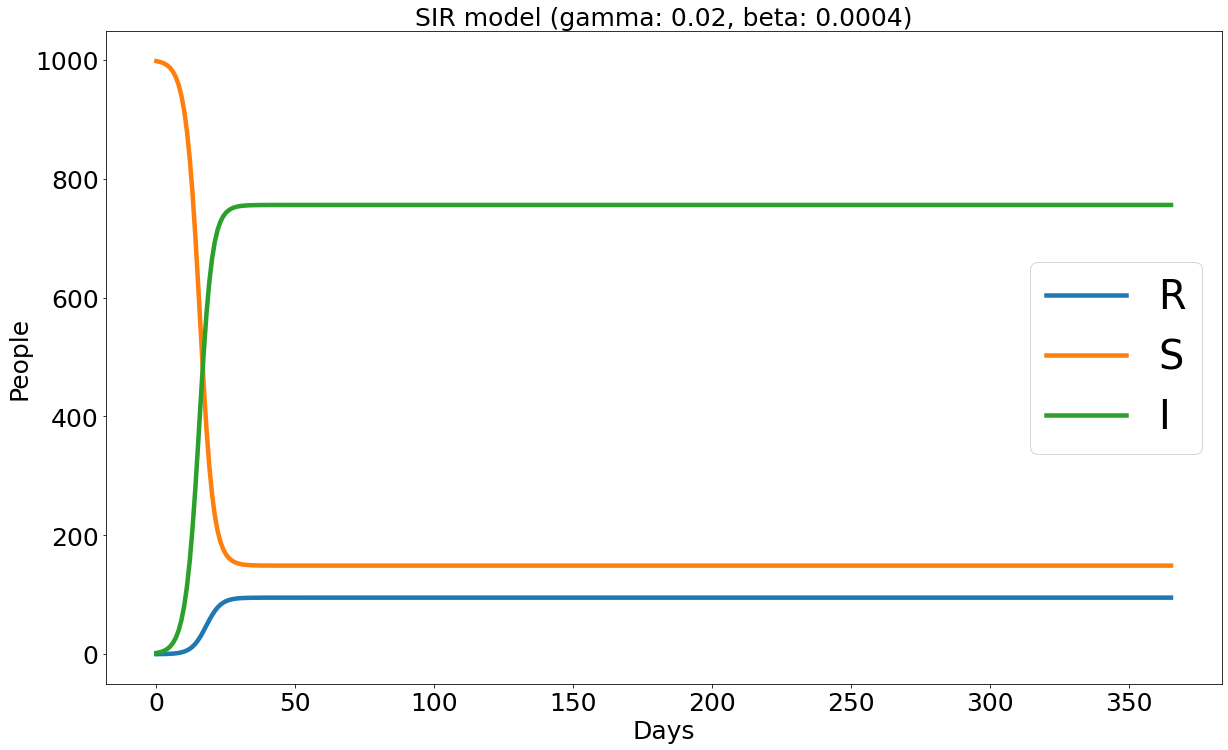}
\par\end{centering}
}
\par\end{centering}
\begin{centering}
\subfloat[Method \#4 ($P=100$)]{\begin{centering}
\includegraphics[scale=0.18]{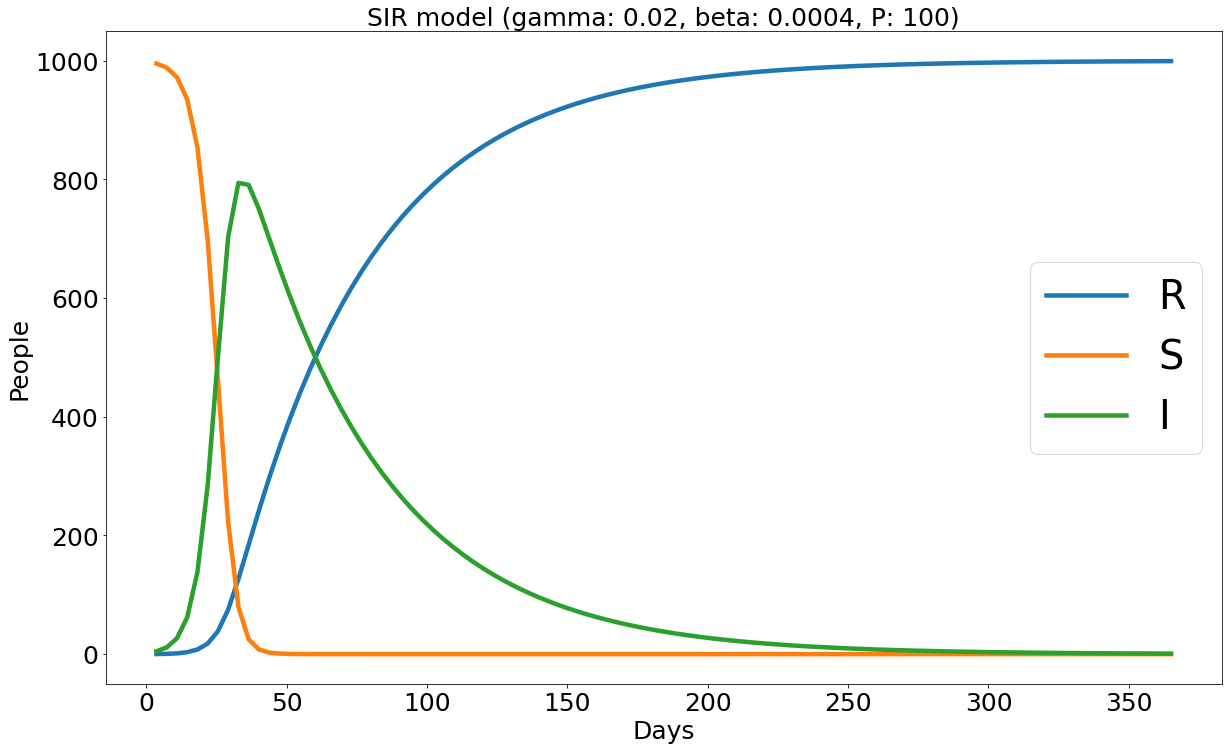}
\par\end{centering}
}\subfloat[Method \#4 ($P=1000$)]{\begin{centering}
\includegraphics[scale=0.18]{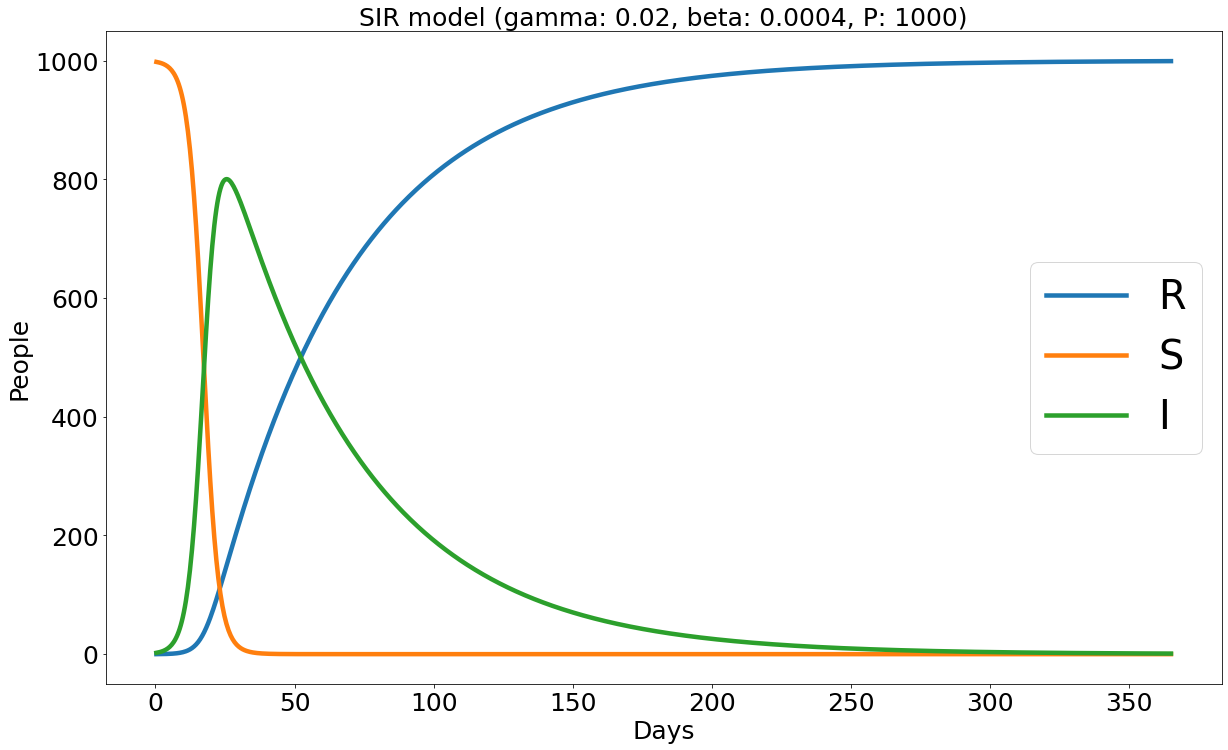}
\par\end{centering}
}
\par\end{centering}
\caption{Graphical illustrations of Test 1. Row 1: Euler-relaxation method.
Row 2: regular linearization method. Row 3: approximate analytic method.
Row 4: direct Euler method.\label{fig:test1}}

\end{figure*}
\end{widetext}

A comparison between methods \#1 and \#2 reveals that while the regular
linearization technique can provide a satisfactory estimate of $I_{\max}$
(800 when considering $P=100$ and $K=2$), method \#2 suffers from
severe numerical instability as illustrated in the second row of Figure
\ref{fig:test1}, particularly when increasing $K$ to obtain a more
accurate graphical representation.

Furthermore, when $P$ and $K$ are relatively small, our proposed
method shows a slight improvement over the conventional Euler method
(method \#4) within the same Euler context. At a coarse grid level,
method \#4 yields a relative error of 0.875\%, while our method achieves
a lower relative error of 0.375\%. Upon increasing $P$ to 1000, both
methods \#1 (with an increased $K=50$) and \#4 demonstrate comparable
accuracy in terms of amplitude and graphical representation, as depicted
in the first and last rows of Figure \ref{fig:test1}.

Our numerical investigation reveals that the true peak value ($I_{\max,\text{true}}$)
is attained on the 24th day by employing sufficiently large values
of $P$ (over 3000) in both reliable methods \#1 and \#4. Comparing
the peak days, it becomes evident from the last row of Table \ref{tab:1}
that our relaxation method outperforms methods \#2 and \#3. While
our method and method \#4 achieve similar accuracy in terms of graphical
simulation and amplitude, our proposed method detects the peak day
earlier and with greater reliability. Specifically, considering small
$P$ and $K$, our relaxation method identifies the peak outbreak
on the 25th day, which closely aligns with the true peak (24th), in
contrast to the peak day of 32nd obtained from method \#4. For larger
$P$, our method predicts an earlier peak occurrence (day 23rd), which
proves advantageous in practical scenarios compared to the peak day
of 25th obtained from method \#4. The ability to predict the peak
event of a disease earlier is of practical significance for decision-makers,
enabling them to implement and sustain timely public health measures
and interventions aimed at mitigating the disease risk.

\subsection*{Test 2}

Our second test focuses on the numerical comparison between two approaches:
applying the well-known Runge-Kutta RK4 method (referred to as method
\#5) to our relaxation scheme (\ref{eq:5}) and applying it directly
to the nonlinear differential equation (\ref{eq:3}) (denoted as method
\#6). Additionally, we compare the convergence speed of method \#5
with method \#1, referred above to as the Euler-relaxation method
(\ref{eq:dis}).

As RK4 is a fourth-order method, we deliberately choose a large population
size of $N=97.47\times10^{6}$ and a transmission rate of $\beta=3\times10^{-9}$.
Assuming the initial infected population is $a=11$, and the removal
rate remains constant at $\gamma=0.05$ throughout the entire six-month
period ($T=180$), we can calculate that the simulated disease reaches
its peak at $I_{\max,\text{true}}=51367769$; cf. (\ref{eq:I_max}).
Moreover, based on numerical observations with a sufficiently large
value of $P$ (>2000), we find that this peak is reached on the 73rd
day.

\begin{widetext}
\begin{figure*}
\begin{centering}
\subfloat[Method \#5 ($P=50,K=20$)]{\begin{centering}
\includegraphics[scale=0.18]{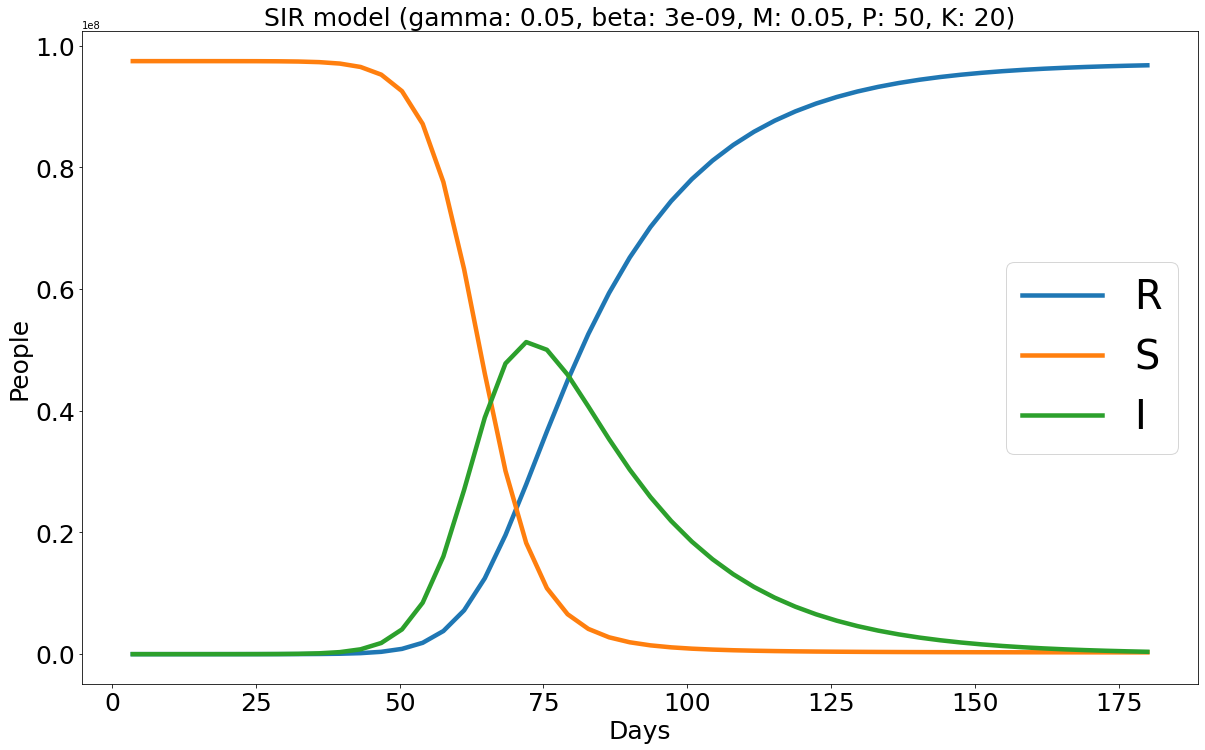}
\par\end{centering}
}\subfloat[Method \#5 ($P=2000,K=50$)]{\begin{centering}
\includegraphics[scale=0.18]{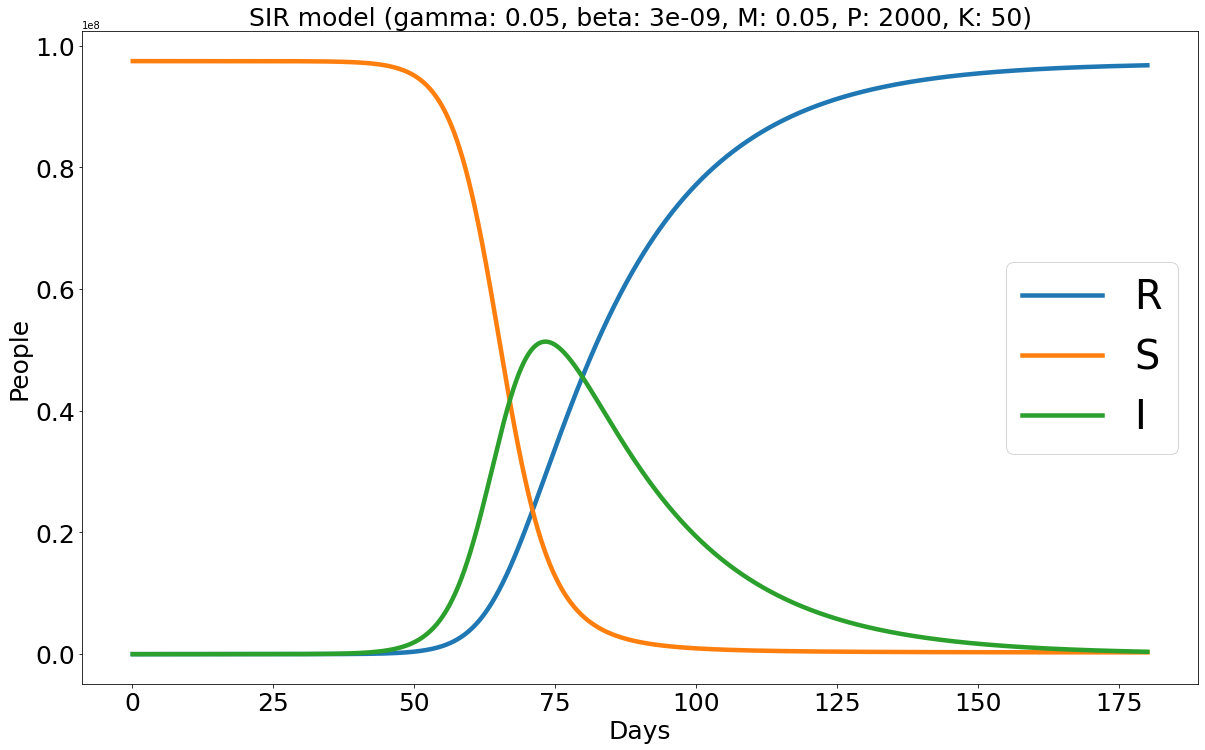}
\par\end{centering}
}
\par\end{centering}
\begin{centering}
\subfloat[Method \#6 ($P=50$)]{\begin{centering}
\includegraphics[scale=0.18]{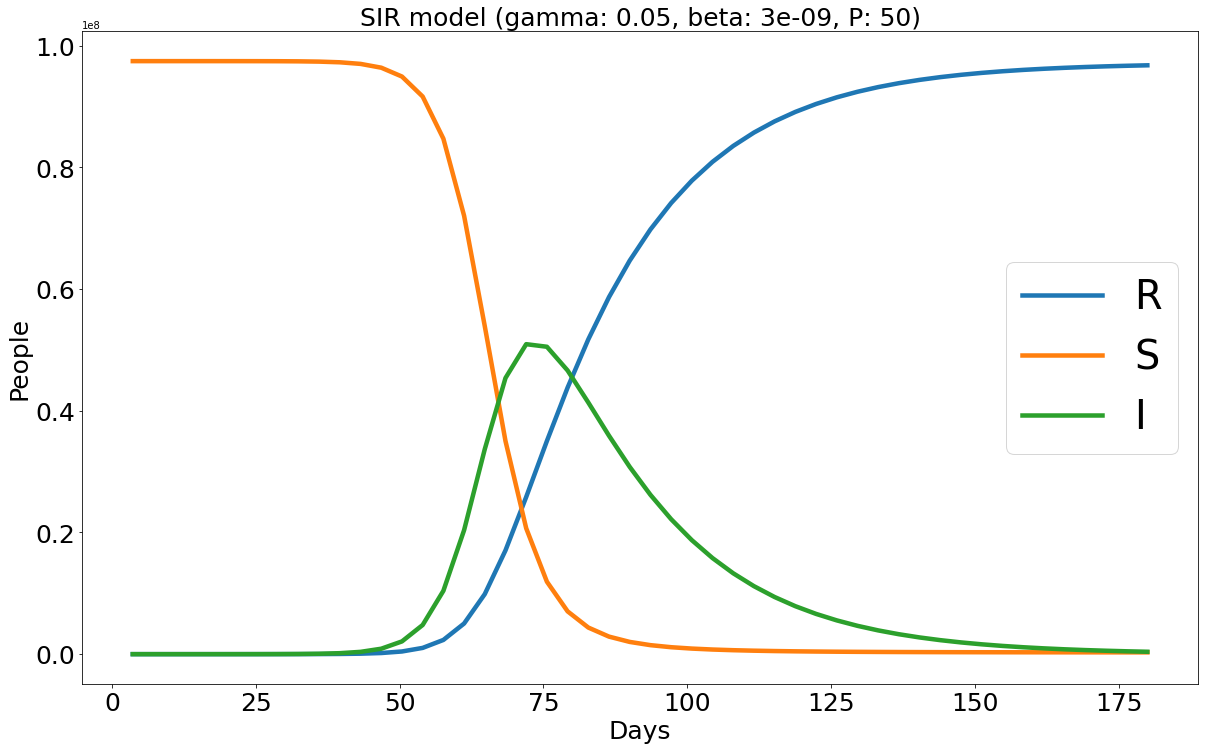}
\par\end{centering}
}\subfloat[Method \#6 ($P=2000$)]{\begin{centering}
\includegraphics[scale=0.18]{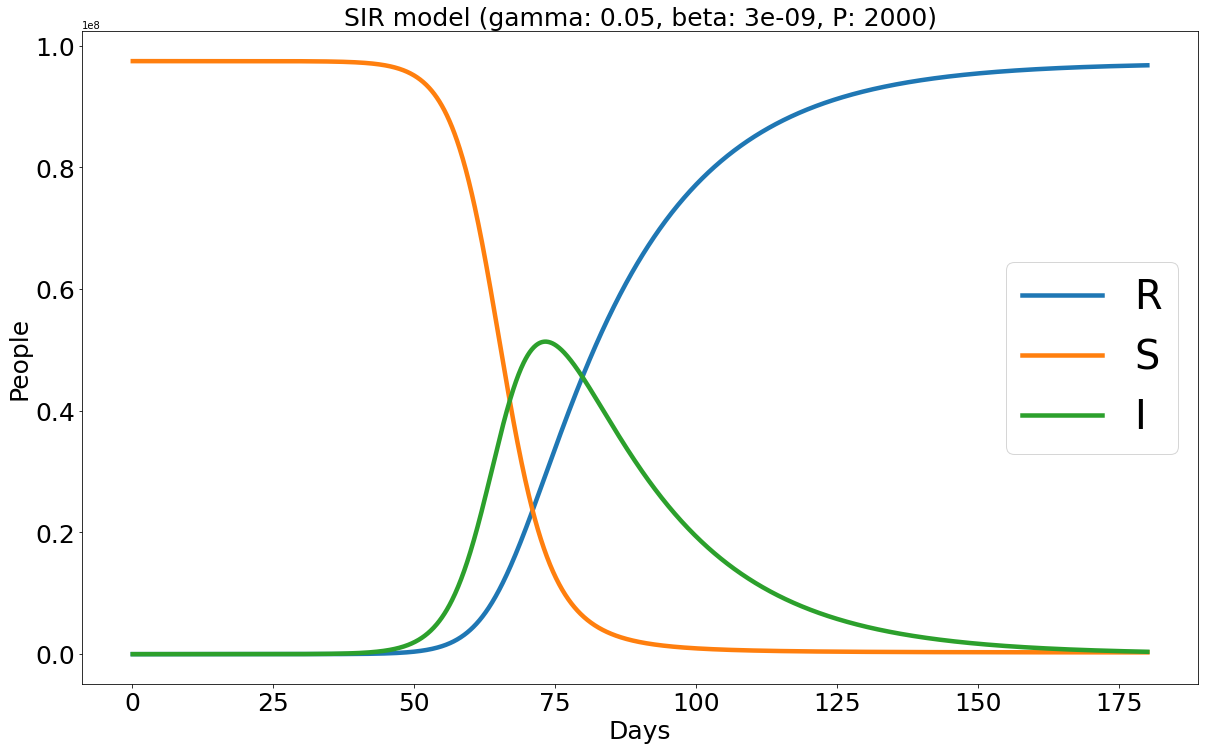}
\par\end{centering}
}
\par\end{centering}
\begin{centering}
\subfloat[Method \#1 ($P=50,K=20$)]{\begin{centering}
\includegraphics[scale=0.18]{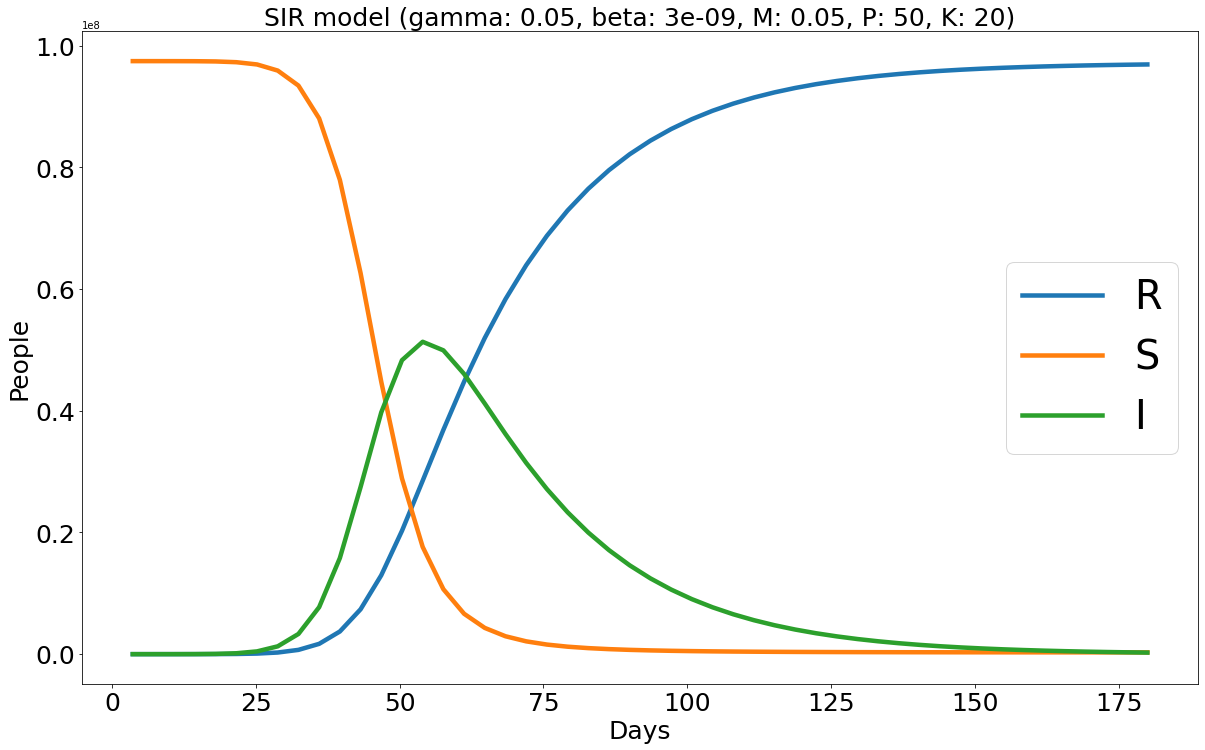}
\par\end{centering}
}\subfloat[Method \#1 ($P=2000,K=50$)]{\begin{centering}
\includegraphics[scale=0.18]{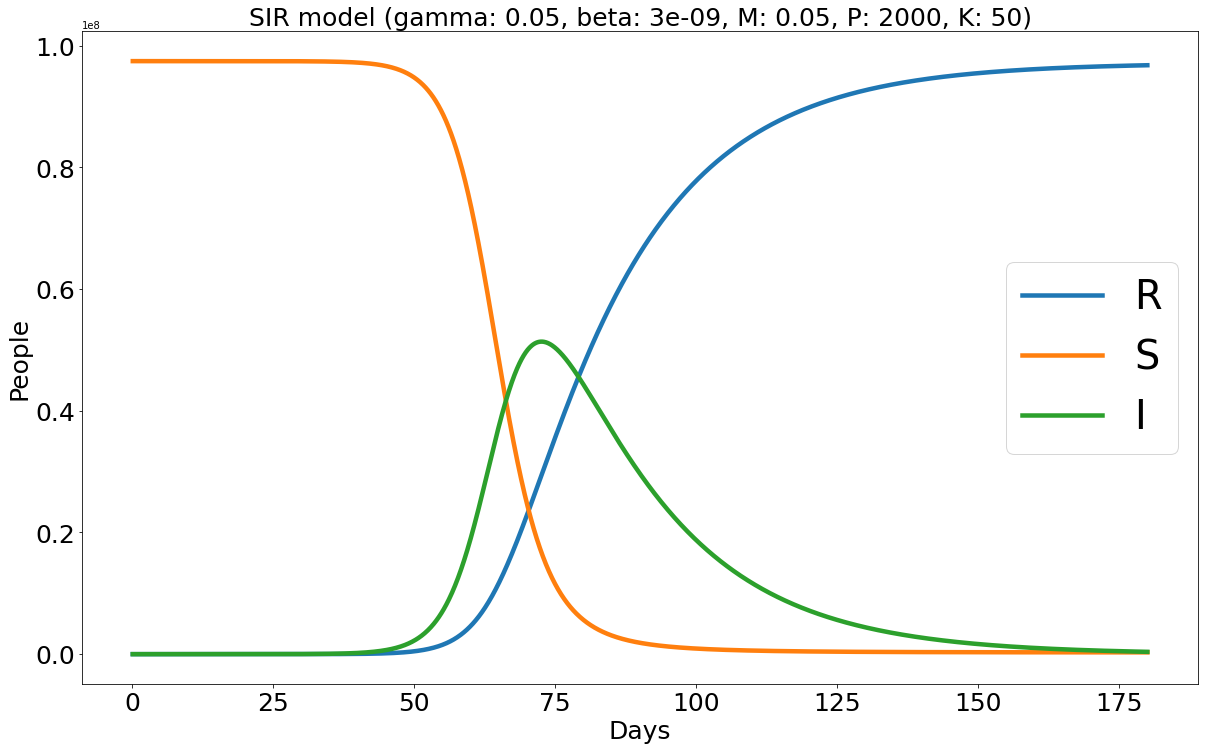}
\par\end{centering}
}
\par\end{centering}
\caption{Graphical illustrations of Test 2. Row 1: RK4-relaxation method. Row
2: direct RK4 method. Row 3: Euler-relaxation method. \label{fig:2}}
\end{figure*}
\end{widetext}

Our numerical results are tabulated in Table \ref{tab:2}, accompanied
by corresponding graphical illustrations in Figure \ref{fig:2}. We
see that within the same RK4 framework, our proposed relaxation method
(method \#5) outperforms the direct approach. When the number of time
steps is small ($P=50$), method \#5 with $K=20$ yields an amplitude
$I_{\max}$ of 51295165 with a relative accuracy of 0.14\%, while
method \#6 achieves 0.81\%. Both methods capture the peak day (72)
well compared to the true value of 73. Note in this case that we choose
$K=20$, a larger value than in Test 1, due to the larger population
under consideration. Cf. Theorem \ref{thm:main}, the choice of $K$
does affect our error estimation which heavily depends on the total
size of the removal population.

We also see that when increasing $P$ to 2000, our proposed method
\#5 with an increased $K=50$ precisely achieves the true amplitude,
$I_{\max,\text{true}}=51367769$, while the direct RK4 method produces
a very close approximation of 51367765. Both methods also identify
the peak day as the 73rd day.

Furthermore, we compare our relaxation method to the Euler and Runge-Kutta
frameworks. In terms of amplitude, although method \#1 initially provides
a better value of 51341234 with an accuracy of 0.05\%, it fails to
accurately detect the peak day, significantly deviating from the true
value of 73 (predicting 54 instead). Increasing $P$ to 2000 improves
the amplitude to 51367573, but it still performs worse than the direct
RK4 method's amplitude of 51367765. Herewith, method \#1 achieves
an improved peak day of 72. Additionally, based on the simulation
of method \#1, we observe that to reach the true amplitude ($I_{\max,\text{true}}=51367769$)
and the true peak day of 73, at least $P=19000$ and $K=100$ are
required. Henceforth, our relaxation method in the RK4 framework,
as readily expected, outperforms itself in the Euler framework.

\begin{table*}[t]
\begin{centering}
\begin{tabular}{|c|c|c|c|c|c|c|}
\hline 
Method \# & 5 & 5 & 6 & 6 & 1 & 1\tabularnewline
\hline 
\# of time step $P$ & $50$ & $2000$ & $50$ & $2000$ & $50$ & $2000$\tabularnewline
\hline 
\# of iteration $K$ & 20 & 50 & None & None & 20 & 50\tabularnewline
\hline 
Amplitude $I_{\max}$ & 51295165 & 51367769 & 50948480 & 51367765 & 51341234 & 51367573\tabularnewline
\hline 
Peak day & 72 & 73 & 72 & 73 & 54 & 72\tabularnewline
\hline 
\end{tabular}
\par\end{centering}
\caption{Values of the computed amplitude $I_{\max}$ obtained from different
methods and the corresponding peak days. Method \#5: our RK4-relaxation
method (\ref{eq:RK4}) applied with (\ref{eq:RR})\textendash (\ref{eq:RRR}).
Method \#6: the conventional RK4 method (\ref{eq:RK4})\textendash (\ref{eq:KK})
applied directly to the nonlinear differential equation (\ref{eq:3}).
Method \#1: our Euler-relaxation method (\ref{eq:dis}). By (\ref{eq:I_max}),
the true amplitude $I_{\max,\text{true}}$ is 51367769 in this scenario.\label{tab:2}}
\end{table*}

\subsection*{Test 3}

As previously mentioned, in our last experiment, we aim to broaden
the scope of the proposed relaxation method by applying it to various
SIR-type models: specifically, the SIRD model and the SIR model with
background mortality. These models share the same input parameters
as those used in Test 1, where we set $N=1000$, $n=998$, $T=365$,
$\gamma=0.02$, $\beta=0.0004$. In the SIRD model (\ref{eq:3.1}),
we choose a death rate of $\sigma=0.01$, which implies a choice of
the relaxation parameter $\overline{M}\ge\gamma+\sigma=0.03$. In
the SIR model with background mortality (\ref{eq:3.13-1}), we use
the background death of $\sigma=0.001$ and select $\hat{M}\ge\gamma=0.02$.

Our numerical findings for the SIRD model are detailed in Figure \ref{fig:3}.
We specifically investigate the scenario where $\overline{M}=0.015$,
thereby contravening the relaxation condition ($\overline{M}\ge0.03$).
Consistent with our theorem concerning non-negativity preservation,
we observe that the relaxed solution with $\overline{M}=0.015$, obtained
from both the Euler and RK4 frameworks, fails to maintain non-negativity
over time. This is evident in the first column of Figure \ref{fig:3}.
When $\overline{M}=0.03$ (a case that adheres to the condition),
we also note that the RK4-relaxation outperforms the Euler-relaxation
approach. Given the $I_{\max,\text{true}}=730$ (as formulated in
(\ref{eq:Imax-SIRD})), and with a coarse mesh of $P=200$ and $K=10$,
we observe that the RK4-relaxation produces an amplitude identical
to the true value. In contrast, the Euler-relaxation yields a value
of 729. It is also worth mentioning that the accurate amplitude is
attained when applying the Euler-relaxation with $P=800$ and $K=10$.
\begin{widetext}
\begin{figure*}
\begin{centering}
\subfloat[Violating Euler-relaxation ($P=100,K=5$, $\overline{M}=0.015$)]{\begin{centering}
\includegraphics[scale=0.18]{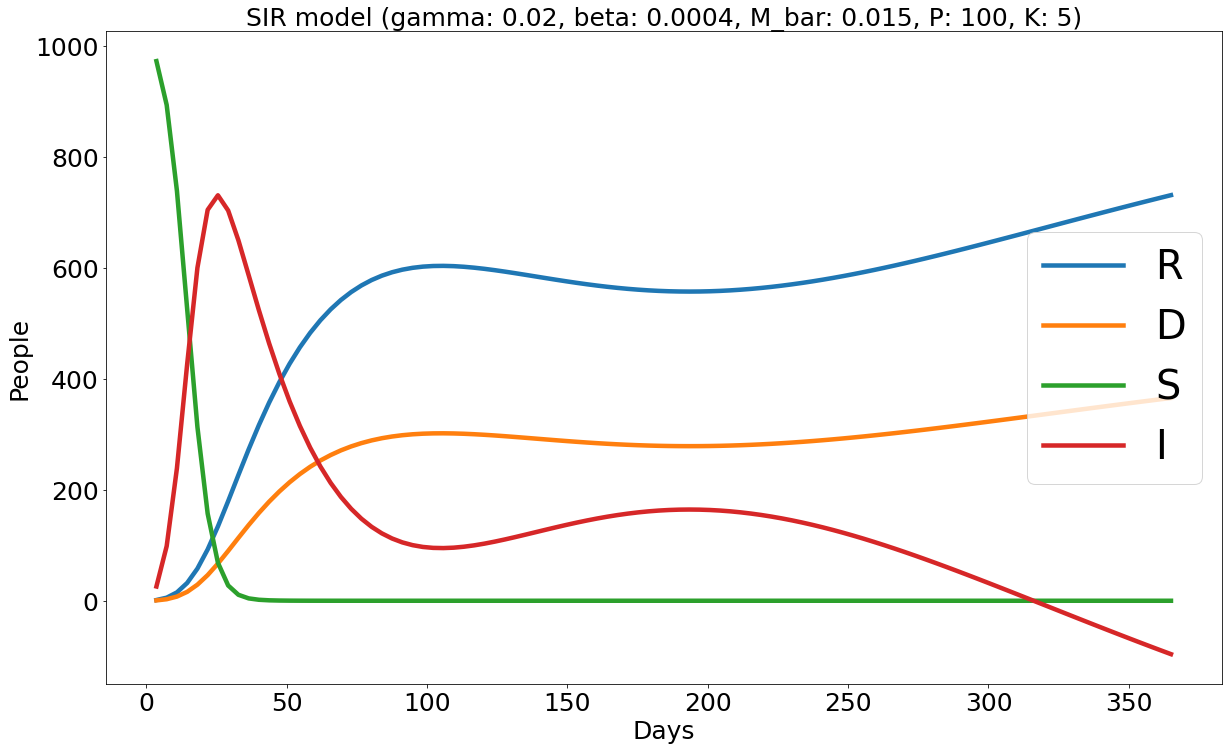}
\par\end{centering}
}\subfloat[Euler-relaxation ($P=200,K=10,\overline{M}=0.03$)]{\begin{centering}
\includegraphics[scale=0.18]{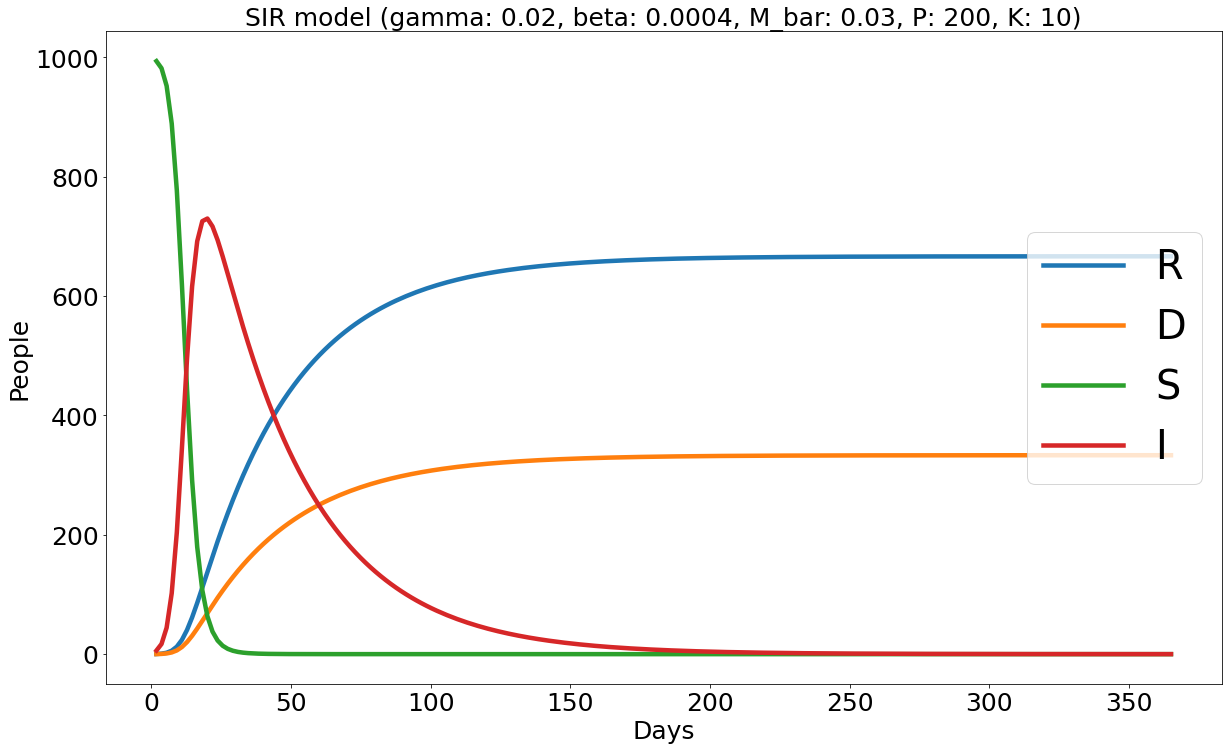}
\par\end{centering}
}
\par\end{centering}
\begin{centering}
\subfloat[Violating RK4-relaxation ($P=100,K=5,\overline{M}=0.015$)]{\begin{centering}
\includegraphics[scale=0.18]{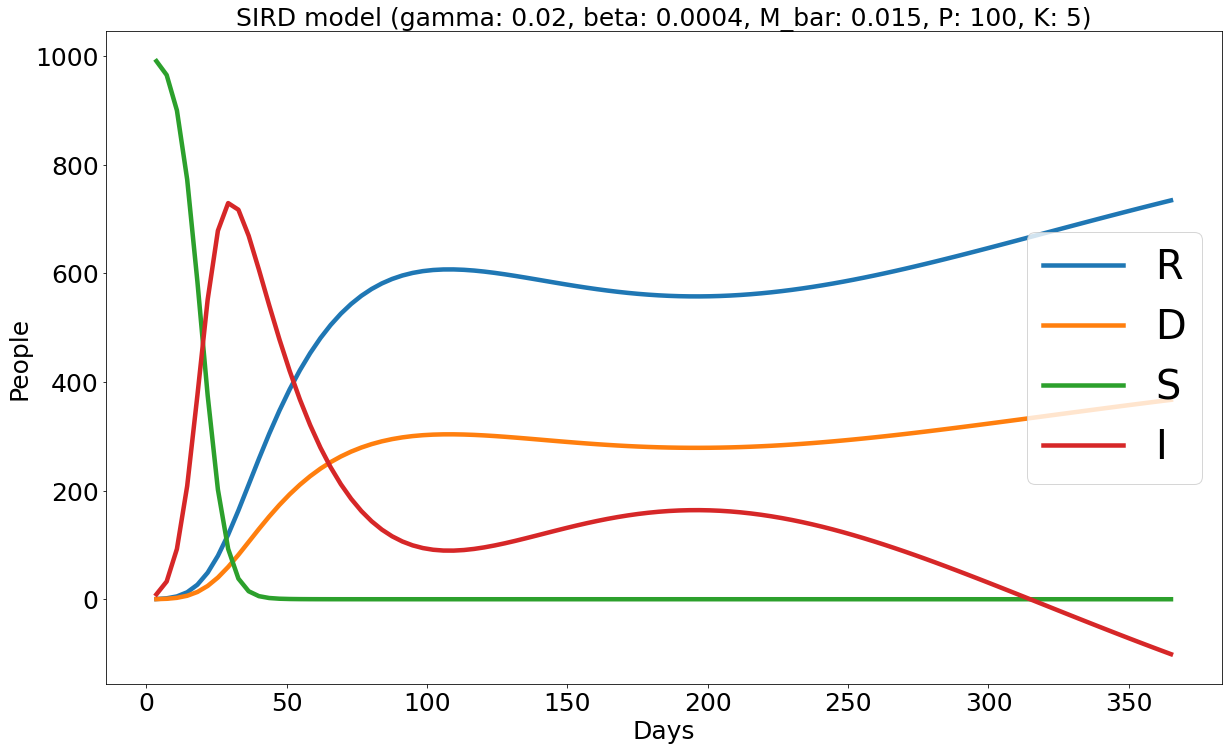}
\par\end{centering}
}\subfloat[RK4-relaxation ($P=200,K=10,\overline{M}=0.03$)]{\begin{centering}
\includegraphics[scale=0.18]{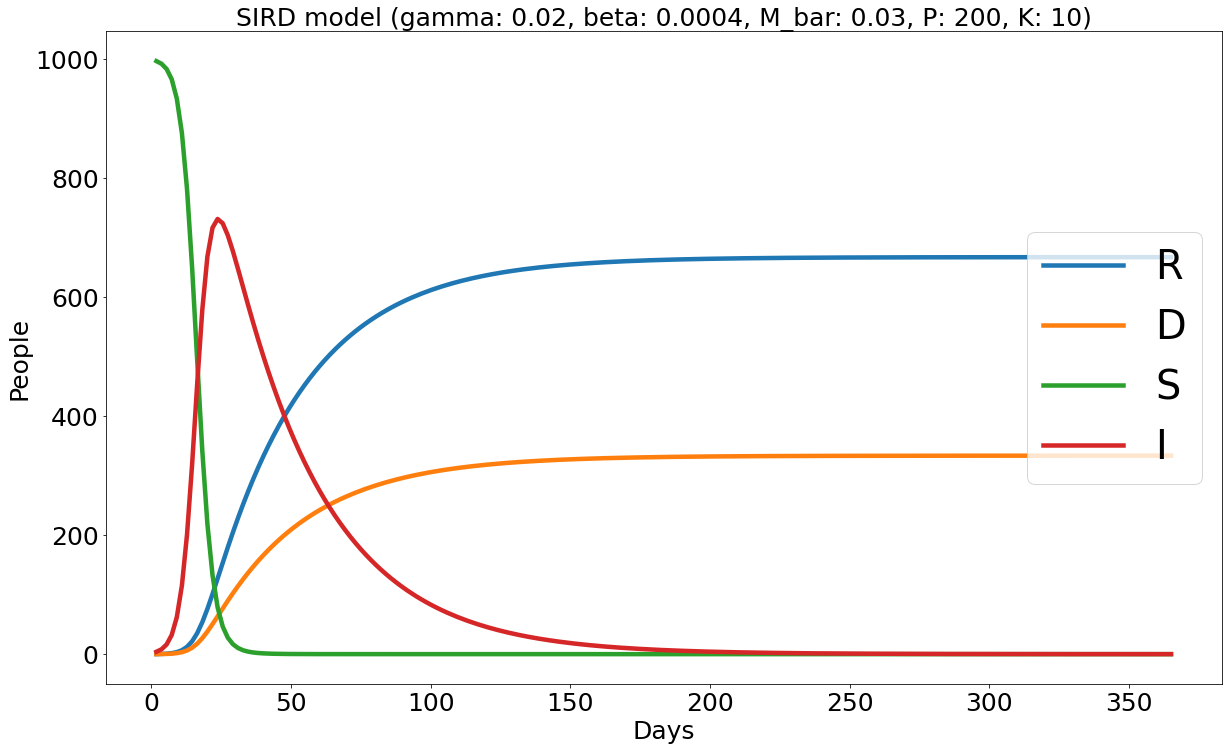}
\par\end{centering}
}
\par\end{centering}
\caption{Graphical illustrations of Test 3 \textendash{} SIRD model. Row 1:
Euler-relaxation method with violating and non-violating cases. Row
2: RK4-relaxation method with violating and non-violating cases. These
illustrations serve to highlight a crucial observation: when the relaxation
parameter is not suitably selected, the numerical solution loses its
adherence to non-negativity preservation.\label{fig:3}}
\end{figure*}

Our numerical results pertaining to the SIR model with background
mortality are presented in Figure \ref{fig:4}. It is evident that
both the Euler and RK4 relaxation methods show numerical stability
and non-negativity preservation as we increase the values of $\left(P,K\right)$
from $\left(100,5\right)$ to $\left(1000,50\right)$.

Consistent with our prior tests, the RK4-relaxation method continues
to outperform the Euler-relaxation method. Leveraging this numerical
stability, we run the RK4-relaxation method using large values of
$P$ and $K$ to determine the numerical amplitude and peak day. Our
findings reveal a numerical amplitude of 777, peaking on the 24th
day. Within the RK4 framework, achieving this numerical amplitude
and peak day requires approximately $P=300$ and $K=20$. In contrast,
the Euler framework demands a minimum of $P=1700$ and $K=20$ for
similar outcomes.

\begin{figure*}
\begin{centering}
\subfloat[Method \#1 ($P=100,K=5$)]{\begin{centering}
\includegraphics[scale=0.18]{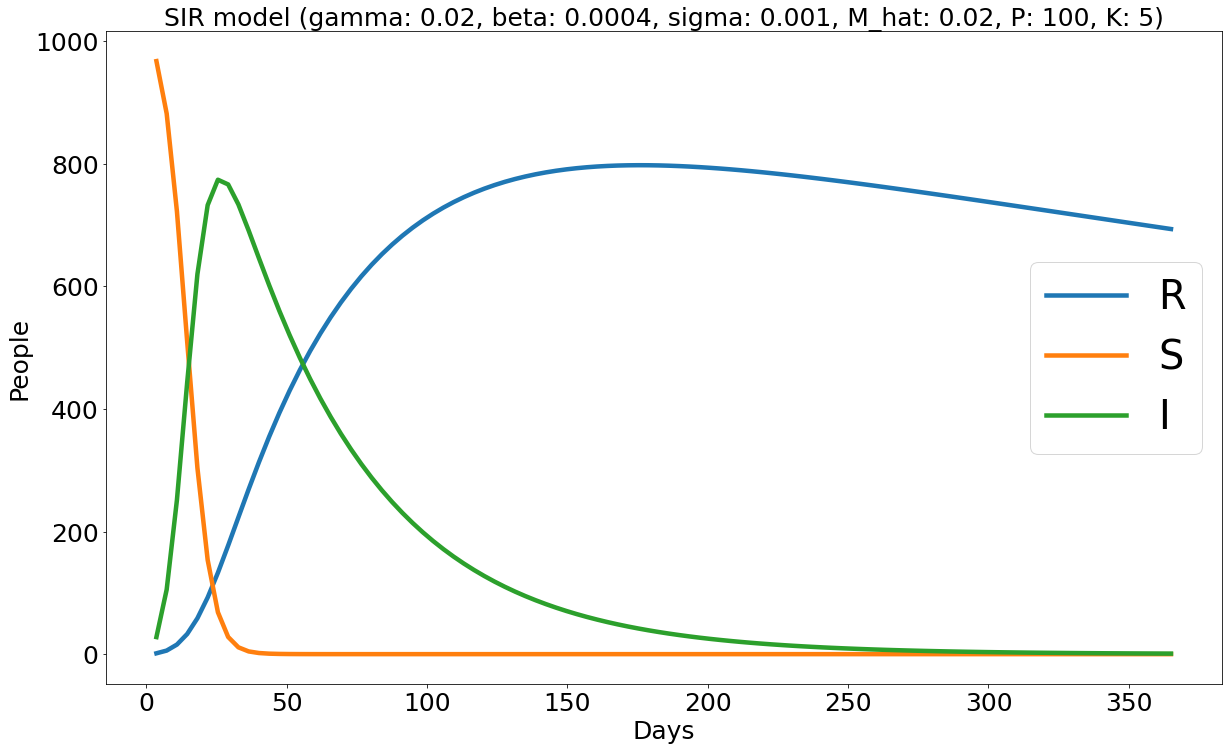}
\par\end{centering}
}\subfloat[Method \#1 ($P=1000,K=50$)]{\begin{centering}
\includegraphics[scale=0.18]{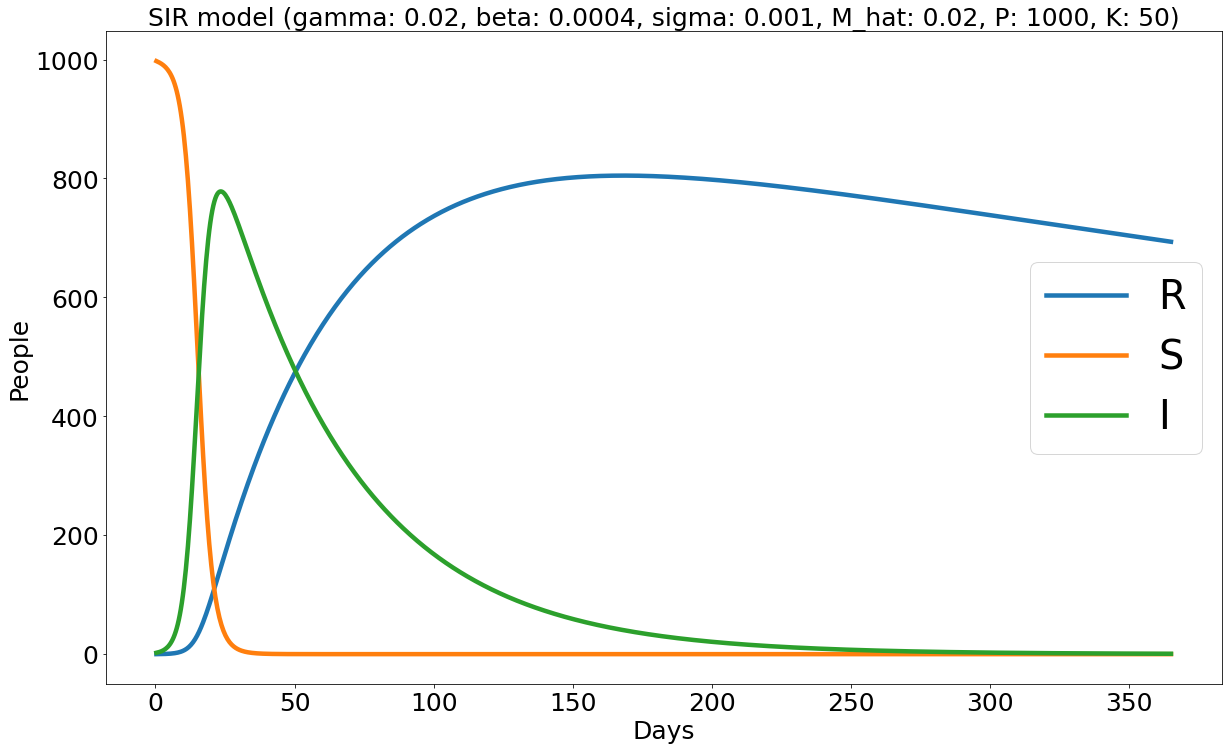}
\par\end{centering}
}
\par\end{centering}
\begin{centering}
\subfloat[Method \#5 ($P=100,K=5$)]{\begin{centering}
\includegraphics[scale=0.18]{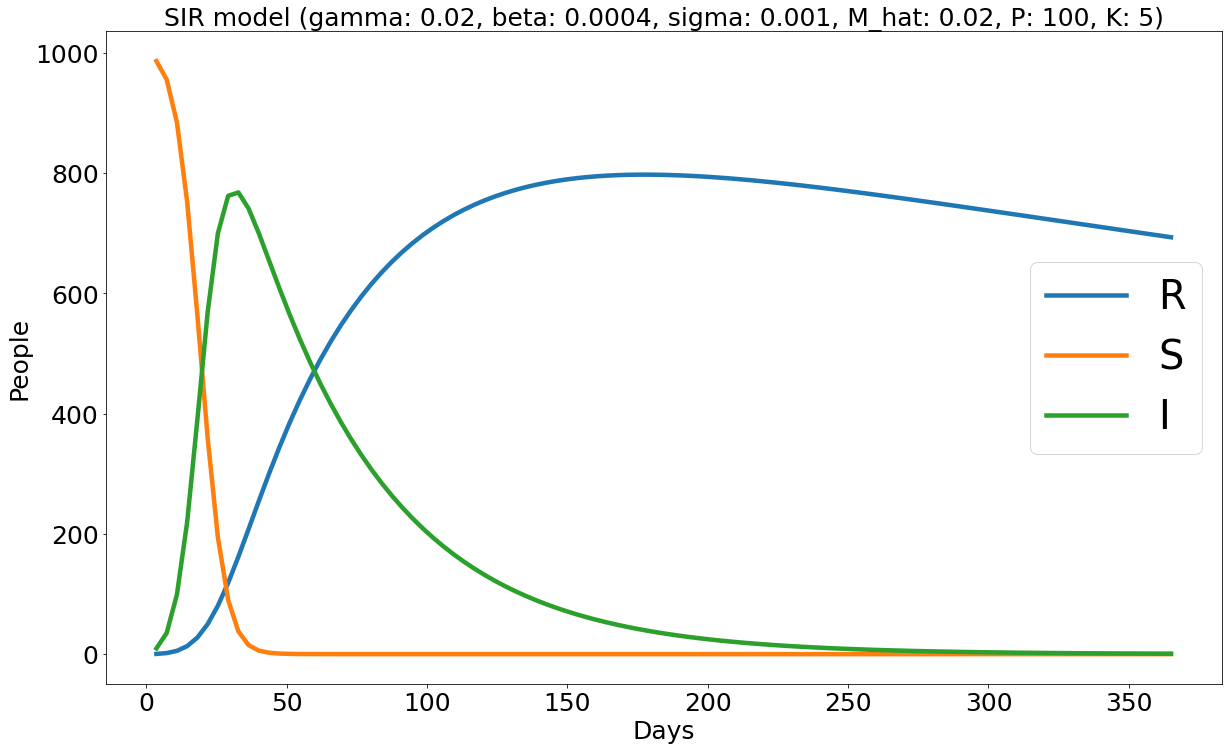}
\par\end{centering}
}\subfloat[Method \#5 ($P=1000,K=50$)]{\begin{centering}
\includegraphics[scale=0.18]{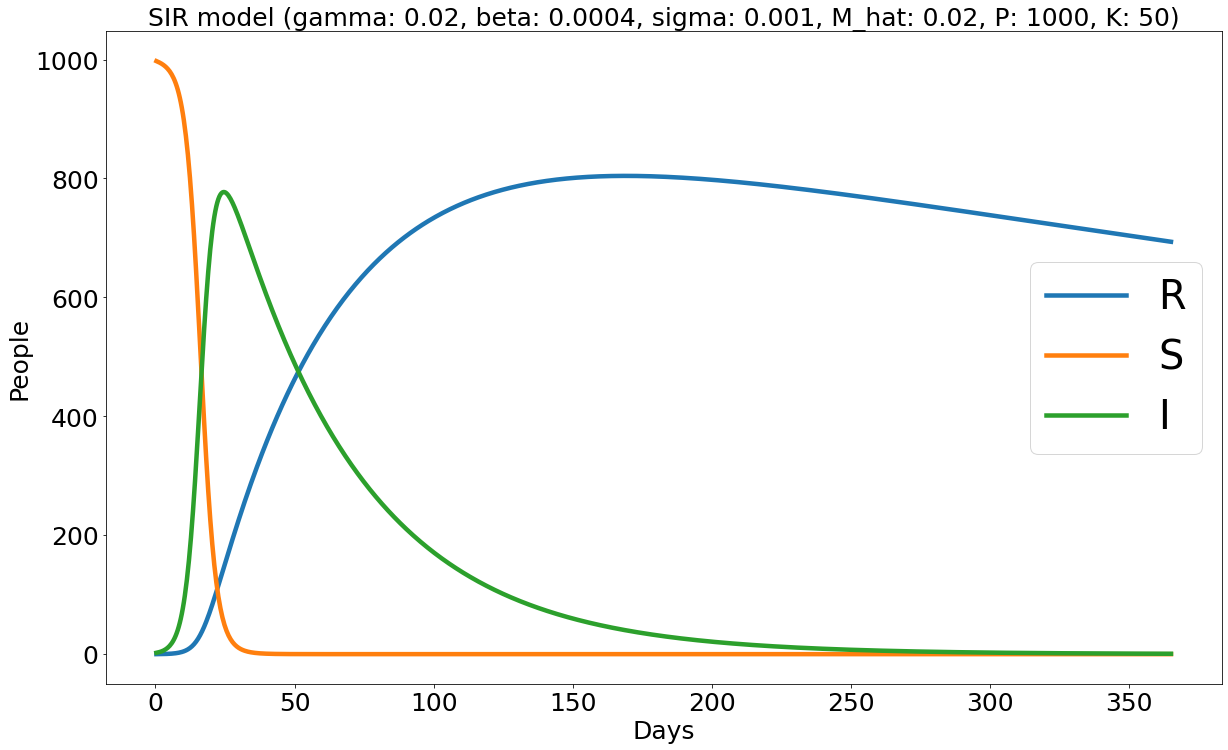}
\par\end{centering}
}
\par\end{centering}
\caption{Graphical illustrations of Test 3 \textendash{} SIR model with background
death. Row 1: Euler-relaxation method with different values of $P$
and $K$. Row 2: RK4-relaxation method with diverse $P$ and $K$
values. In these illustrations, numerical stability and non-negativity
preservation are observed.\label{fig:4}}
\end{figure*}
\end{widetext}

\section{Concluding remarks\label{sec:45}}

This work presents a novel numerical approach for solving the SIR
model in population dynamics. While various approximation methods
have been proposed for this classical model, the analysis of their
convergence has been limited and challenging. Our approach introduces
the relaxation procedure to approximate the continuous model. By carefully
selecting the relaxation parameter, we achieve global strong convergence
of the scheme and effectively preserve non-negativity. The proposed
scheme is explicit and straightforward to implement, enabling us to
obtain the approximate solution at either the discrete or analytical
level. Additionally, we showcase the applicability of our scheme to
numerous variants of the SIR model.


In our future work, we aim to extend our method to a globally strongly
convergent higher-order scheme. Furthermore, we plan to apply the
method to more complex SIR-based models involving multiple compartments.

\section*{Credit authorship contribution statement}

\textbf{V. A. Khoa:} Supervision, Conceptualization, Formal analysis,
Writing\textendash review \& editing. \textbf{P. M. Quan:} Software,
Visualization. \textbf{K. W. Blayneh:} Formal analysis, Writing\textendash review
\& editing. \textbf{J. Allen:} Formal analysis, Visualization.

\section*{Declaration of competing interest}

The authors declare that they have no known competing financial interests
or personal relationships that could have appeared to influence the
work reported in this paper.

\section*{Data availability}

Simulated data will be made available on request.

\section*{Acknowledgments}

This research has received funding from the National Science Foundation
(NSF). Specifically, V. A. K., P. M. Q., and J. A. extend their gratitude
for the invaluable support provided by NSF. V. A. K. and J. A. also
hold deep appreciation for the Florida A\&M University Rattler Research
program, its esteemed committee, and Dr. Tiffany W. Ardley. Their
unwavering dedication has facilitated an exceptional academic journey
for the mentee (J. A.) and the mentor (V. A. K.).

Furthermore, V. A. K. wishes to express many thanks to Dr. Charles
Weatherford (Florida A\&M University) and Dr. Ziad Musslimani (Florida
State University). Their support has been instrumental in shaping
V. A. K.'s early research career. Lastly, this work is complete on
a momentous personal milestone \textendash{} the wedding day of V.
A. K. and the bride, Huynh Thi Kim Ngan.

\section*{Appendix}

\subsection*{Proof of Theorem \ref{thm:main}}
\begin{widetext}
\uline{Step 1:} Define $\mathcal{E}_{k}\left(t\right)=R_{k}\left(t\right)-R\left(t\right)$
for $k=1,2,3,\ldots$. It follows from (\ref{eq:5}) and (\ref{eq:3})
that $\mathcal{E}_{k}$ satisfies the following differential equation:

\begin{align}
\mathcal{E}_{k}'\left(t\right)+M\mathcal{E}_{k}\left(t\right) & =-g\left(R_{k-1}\left(t\right)\right)+g\left(R\left(t\right)\right)+M\mathcal{E}_{k-1}\left(t\right)\nonumber \\
 & =p\left(R\left(t\right)\right)-p\left(R_{k-1}\left(t\right)\right),\label{eq:6}
\end{align}
where we have denoted $p\left(r\right)=g\left(r\right)-Mr$ for $r\ge0$.
Herewith, by Theorems \ref{thm:1} and \ref{thm:2}, we are allowed
to consider $r\ge0$. We can compute that $p'\left(r\right)=-\gamma n\mu e^{-\mu r}+\gamma-M$.
Then, for $M\ge\gamma$ and since $r\ge0$, we estimate that
\[
\gamma-\gamma n\mu-M\le p'\left(r\right)=-\gamma n\mu e^{-\mu r}+\gamma-M<0,
\]
which shows
\begin{equation}
\left|p'\left(r\right)\right|\le M+\gamma n\mu-\gamma.\label{eq:7}
\end{equation}
Therefore, the left-hand side of (\ref{eq:6}) can be bounded from
above by
\[
\mathcal{E}_{k}'\left(t\right)+M\mathcal{E}_{k}\left(t\right)\le\left(M+\gamma n\mu-\gamma\right)\left|\mathcal{E}_{k-1}\left(t\right)\right|.
\]
Using the H\"older inequality, we find that
\begin{align*}
e^{2Mt}\left|\mathcal{E}_{k}\left(t\right)\right|^{2} & \le\left(M+\gamma n\mu-\gamma\right)^{2}\left(\int_{0}^{t}e^{Ms}\left|\mathcal{E}_{k-1}\left(s\right)\right|ds\right)^{2}\\
 & \le\left(M+\gamma n\mu-\gamma\right)^{2}\int_{0}^{t}e^{2Ms}ds\int_{0}^{t}\left|\mathcal{E}_{k-1}\left(s\right)\right|^{2}ds.
\end{align*}
Thus, we deduce that
\[
\left|\mathcal{E}_{k}\left(t\right)\right|^{2}\le\frac{1}{2M}\left(M+\gamma n\mu-\gamma\right)^{2}\left(1-e^{-2Mt}\right)\int_{0}^{t}\left|\mathcal{E}_{k-1}\left(s\right)\right|^{2}ds.
\]
By the elementary inequality $e^{-x}+x\ge1$, we obtain the following
estimate
\begin{equation}
\left|\mathcal{E}_{k}\left(t\right)\right|^{2}\le\left(M+\gamma n\mu-\gamma\right)^{2}t\int_{0}^{t}\left|\mathcal{E}_{k-1}\left(s\right)\right|^{2}ds.\label{eq:2.77}
\end{equation}
\uline{Step 2:} By induction, we can show that for any $2\le k\in\mathbb{N}$
\begin{equation}
\left|\mathcal{E}_{k}\left(t\right)\right|^{2}\le\left(M+\gamma n\mu-\gamma\right)^{2k}t\int_{0}^{t}s_{1}\int_{0}^{s_{1}}\ldots s_{k-1}\int_{0}^{s_{k-1}}\left|\mathcal{E}_{0}\left(s_{k}\right)\right|^{2}ds_{k}ds_{k-1}\ldots ds_{1}\label{eq:2.7}
\end{equation}
It follows from (\ref{eq:2.77}) that (\ref{eq:2.7}) holds true for
$k=2$. Indeed,
\[
\left|\mathcal{E}_{2}\left(t\right)\right|^{2}\le\left(M+\gamma n\mu-\gamma\right)^{2}t\int_{0}^{t}\left|\mathcal{E}_{1}\left(s_{1}\right)\right|^{2}ds_{1}\le\left(M+\gamma n\mu-\gamma\right)^{4}t\int_{0}^{t}s_{1}\int_{0}^{s_{1}}\left|\mathcal{E}_{0}\left(s_{2}\right)\right|^{2}ds_{2}ds_{1}.
\]
Assume that (\ref{eq:2.7}) holds true for $k=k_{0}$. We show that
it also holds true for $k=k_{0}+1$. By (\ref{eq:2.77}), we have
\begin{align*}
 & \left|\mathcal{E}_{k_{0}+1}\left(t\right)\right|^{2}\\
 & \le\left(M+\gamma n\mu-\gamma\right)^{2}t\int_{0}^{t}\left|\mathcal{E}_{k_{0}}\left(s\right)\right|^{2}ds\\
 & \le\left(M+\gamma n\mu-\gamma\right)^{2}t\int_{0}^{t}\left(M+\gamma n\mu-\gamma\right)^{2k_{0}}s\int_{0}^{s}s_{1}\int_{0}^{s_{1}}\ldots s_{k_{0}-1}\int_{0}^{s_{k_{0}-1}}\left|\mathcal{E}_{0}\left(s_{k_{0}}\right)\right|^{2}ds_{k_{0}}ds_{k_{0}-1}\ldots ds_{1}ds\\
 & \le\left(M+\gamma n\mu-\gamma\right)^{2\left(k_{0}+1\right)}\int_{0}^{t}s_{1}\int_{0}^{s_{1}}\ldots s_{k_{0}}\int_{0}^{s_{k_{0}}}\left|\mathcal{E}_{0}\left(s_{k_{0}+1}\right)\right|^{2}ds_{k_{0}+1}ds_{k_{0}}\ldots ds_{1}.
\end{align*}
Hence, we complete Step 2.

\uline{Step 3:} By (\ref{eq:2.7}), observe that $0\le s_{k}\le s_{k-1}\le\ldots\le s_{1}\le t$.
Combining this, (\ref{eq:2.7}) and the fact that $R\in C^{1}$ gives
\begin{align*}
\left|\mathcal{E}_{k}\left(t\right)\right|^{2} & \le\left(M+\gamma n\mu-\gamma\right)^{2k}t^{k+1}\max_{0\le t\le T}\left|\mathcal{E}_{0}\left(t\right)\right|^{2}\int_{0}^{t}s_{1}\int_{0}^{s_{1}}\ldots s_{k-1}\int_{0}^{s_{k-1}}ds_{k}ds_{k-1}\ldots ds_{1}\\
 & \le\left(M+\gamma n\mu-\gamma\right)^{2k}\frac{t^{k+1}}{k!}\max_{0\le t\le T}\left|\mathcal{E}_{0}\left(t\right)\right|^{2}.
\end{align*}
Note that we have the $k$ and time independence of $M+\gamma n\mu-\gamma$
and $t\le T$. Moreover, we know that $\mathcal{E}_{0}\left(t\right)=R_{0}\left(t\right)-R\left(t\right)=-R\left(t\right)$
by the choice $R_{0}\left(t\right)=0$. Therefore, in view of the
fact that $\lim_{k\to\infty}\frac{Q^{k}}{k!}=0$ for any constant
$Q>0$, we can always find $\overline{k}>0$ such that for any $k\ge\overline{k}$,
\[
\left(M+\gamma n\mu-\gamma\right)^{2k}\frac{T^{k+1}}{k!}<1.
\]
Hence, we obtain the strong convergence of the sequence $\left\{ R_{k}\right\} _{k=0}^{\infty}$
toward the true solution $R$.
\end{widetext}

\subsection*{Proof of Corollary \ref{cor:4}}
\begin{widetext}
We define $\mathcal{E}_{k}\left(t\right)=R_{k}\left(t\right)-R\left(t\right)$
and $p\left(r\right)=g\left(r\right)-Mr$ as above. Multiplying (\ref{eq:6})
by $\mathcal{E}_{k}\left(t\right)$ and using (\ref{eq:7}) yield
\[
\frac{1}{2}\frac{d}{dt}\mathcal{E}_{k}^{2}\left(t\right)+M\mathcal{E}_{k}^{2}\left(t\right)=\left[p\left(R\left(t\right)\right)-p\left(R_{k-1}\left(t\right)\right)\right]\mathcal{E}_{k}\left(t\right)\le\frac{M+\gamma n\mu-\gamma}{2}\mathcal{E}_{k-1}^{2}\left(t\right)+\frac{M+\gamma n\mu-\gamma}{2}\mathcal{E}_{k}^{2}\left(t\right).
\]
Equivalently, we obtain
\[
\frac{d}{dt}\mathcal{E}_{k}^{2}\left(t\right)+\left(M-\gamma n\mu+\gamma\right)\mathcal{E}_{k}^{2}\left(t\right)\le\left(M+\gamma n\mu-\gamma\right)\mathcal{E}_{k-1}^{2}\left(t\right).
\]
Notice that by the choice $M\ge\gamma$, it holds true that $M>\gamma n\mu-\gamma$
when $n\mu<1$. Using the integrating factor $e^{\left(M-\gamma n\mu+\gamma\right)t}$
and taking integration with respect to $t$, we get
\begin{align*}
\mathcal{E}_{k}^{2}\left(t\right) & \le e^{-\left(M-\gamma n\mu+\gamma\right)t}\left(M+\gamma n\mu-\gamma\right)\int_{0}^{t}e^{\left(M-\gamma n\mu+\gamma\right)s}\mathcal{E}_{k-1}^{2}\left(s\right)ds\\
 & \le e^{-\left(M-\gamma n\mu+\gamma\right)t}\left[e^{\left(M-\gamma n\mu+\gamma\right)t}-1\right]\frac{M+\gamma n\mu-\gamma}{M-\gamma n\mu+\gamma}\max_{0\le t\le T}\mathcal{E}_{k-1}^{2}\left(t\right).
\end{align*}
Henceforth, we obtain
\begin{equation}
\max_{0\le t\le T}\left|\mathcal{E}_{k}\left(t\right)\right|\le\left(\frac{M+\gamma n\mu-\gamma}{M-\gamma n\mu+\gamma}\right)^{1/2}\max_{0\le t\le T}\left|\mathcal{E}_{k-1}\left(t\right)\right|.\label{eq:8}
\end{equation}
By induction and the fact that $R_{0}\left(t\right)=0$, we deduce
\[
\max_{0\le t\le T}\left|\mathcal{E}_{k}\left(t\right)\right|\le\left(\frac{M+\gamma n\mu-\gamma}{M-\gamma n\mu+\gamma}\right)^{k/2}\max_{0\le t\le T}\left|\mathcal{E}_{0}\left(t\right)\right|=\left(\frac{M+\gamma n\mu-\gamma}{M-\gamma n\mu+\gamma}\right)^{k/2}\max_{0\le t\le T}\left|R\left(t\right)\right|.
\]
Since $M+\gamma n\mu-\gamma<M-\gamma n\mu+\gamma$ when $n\mu<1$,
we obtain the target estimate (\ref{eq:2.9}).
\end{widetext}

\bibliographystyle{plain}
\bibliography{mybib}

\end{document}